\documentclass[11pt, reqno, psamsfonts]{amsart}
\pdfoutput=1

\usepackage{amssymb}
\usepackage{amsthm}
\usepackage{amsmath}
\usepackage{latexsym}
\usepackage[T1]{fontenc}
\usepackage[utf8]{inputenc}
\usepackage[russian, french, english]{babel}

\usepackage{graphicx}
\usepackage{wrapfig}
\usepackage[justification=centering, labelfont=bf]{caption}
\usepackage{mathtools}
\usepackage{tikz}
\usepackage{amsbsy}
\usepackage[inline]{enumitem}
\usepackage{mathrsfs}
\usetikzlibrary{shapes,snakes}
\usetikzlibrary{arrows.meta}
\usetikzlibrary{decorations.pathmorphing}
\usetikzlibrary{patterns}
\usepackage{float}
\usepackage{array}
\usepackage{multicol}
\usepackage{stmaryrd}
\usepackage{eucal,eufrak}

\usepackage{caladea}
\usepackage[T1]{fontenc}
\usepackage[euler-digits,small]{eulervm}
\usepackage[bb=pazo,scr=esstix]{mathalfa}

\usepackage{upgreek}
\usepackage{titlesec}
\usepackage[spacing=true,kerning=true,babel=true,tracking=true]{microtype}
\usepackage[shortcuts]{extdash}
\usepackage[foot]{amsaddr} 
\usepackage{framed}
\usepackage{algpseudocode}
\usepackage{algorithm}
\definecolor{shadecolor}{gray}{0.95}
\usepackage[left=1in,right=1in,top=1in,bottom=1in,bindingoffset=0cm]{geometry}
\usepackage{xspace}
\usepackage{bm}
\usepackage{centernot}

\usepackage[backend=biber, style=alphabetic, sorting=nyt, maxnames=100,backref=true]{biblatex}
\usepackage[hidelinks]{hyperref}
\title{\lsstyle Probabilistic constructions in continuous combinatorics \\ and a bridge to distributed algorithms}
\date{}
\author{\lsstyle Anton~Bernshteyn}
\address{\textls{\normalfont{}School of Mathematics, Georgia Institute of Technology, Atlanta, GA, USA}}
\email{bahtoh@gatech.edu}
\thanks{This research was partially supported by the NSF grant DMS-2045412.}

\newtheoremstyle{bfnote}%
{}{}%
{\slshape}{}%
{\bfseries}{\bfseries.}%
{ }%
{\thmname{#1}\thmnumber{ #2}\thmnote{ \ep{\normalfont{}#3}}}

\newtheoremstyle{defbfnote}%
{}{}%
{}{}%
{\bfseries}{.}%
{ }%
{\thmname{#1}\thmnumber{ #2}\thmnote{ (#3)}}

\newtheoremstyle{claim}%
{}{}%
{\slshape}{}%
{\itshape}{.}%
{ }%
{\thmname{#1}\thmnumber{ #2}\thmnote{ \ep{\normalfont{}#3}}}

\newtheoremstyle{smalldefn}%
{}{}%
{}{}%
{\itshape}{.}%
{ }%
{\thmname{#1}\thmnumber{ #2}\thmnote{ \ep{\normalfont{}#3}}}

\theoremstyle{bfnote}

\newtheorem{theo}[equation]{Theorem}
\newtheorem*{theo*}{Theorem}

\newtheorem{lemma}[equation]{Lemma}
\newtheorem{corl}[equation]{Corollary}

\newtheorem*{claim*}{Claim}
\newtheorem{big_claim}[equation]{Claim}
\newtheorem*{corl*}{Corollary}

\theoremstyle{claim}

\theoremstyle{smalldefn}

\newenvironment{stepproof}{\noindent$\rhd$\hspace{1em}}{\hfill$\blacktriangleleft$\smallskip}

\newcommand*{\myproofname}{Proof}

\theoremstyle{definition}
\newtheorem{defn}[equation]{Definition}
\newtheorem*{defn*}{Definition}

\newtheorem{ques}[equation]{Question}

\newtheorem*{exmp*}{Example}

\newtheorem*{assum*}{Assumptions}

\theoremstyle{remark}
\newtheorem*{ques*}{Question}
\newtheorem*{remk*}{Remark}

\newcommand{\0}{\varnothing}
\newcommand{\set}[1]{\{#1\}}

\newcommand{\dom}{\mathrm{dom}}

\newcommand{\acts}{\curvearrowright}
\newcommand{\N}{{\mathbb{N}}}
\newcommand{\Z}{\mathbb{Z}}

\renewcommand{\epsilon}{\varepsilon}
\renewcommand{\phi}{\varphi}
\renewcommand{\theta}{\vartheta}
\renewcommand{\leq}{\leqslant}
\renewcommand{\geq}{\geqslant}

\newcommand{\fins}[1]{[#1]^{<\infty}}

\renewcommand{\G}{\Gamma}

\newcommand{\defeq}{\coloneqq}
\newcommand{\rest}[2]{{{#1}\vert_{#2}}}
\newcommand{\emphd}[1]{{\fontseries{b}\selectfont{#1}}}
\newcommand{\Stab}{\mathrm{St}}
\renewcommand{\P}{\mathbb{P}}

\newcommand{\Nbhd}{{N}}

\newcommand{\LOCAL}{$\mathsf{LOCAL}$\xspace}
\newcommand{\B}{\mathscr{B}}
\newcommand{\Pat}{\mathscr{P}}

\newcommand{\pr}{\mathsf{p}}
\newcommand{\de}{\mathsf{d}}
\newcommand{\vdeg}{\mathsf{vdeg}}
\newcommand{\ord}{\mathsf{ord}}

\newcommand{\Free}{\mathrm{Free}}

\newcommand{\bemph}[1]{{\normalfont#1}} 
\newcommand{\ep}[1]{\bemph{(}#1\bemph{)}} 

\newenvironment{scproof}[1][Proof]{\begin{proof}[\lsstyle\itshape#1]}{\end{proof}}

\numberwithin{equation}{section}

\renewcommand{\thesubsection}{\arabic{section}.\Alph{subsection}}

\usepackage{etoolbox}

\titleformat{\section}[block]{\scshape\filcenter}{\thesection.}{1ex}{}
\titleformat{\subsection}[block]{\bfseries\filcenter}{\thesubsection.}{1ex}{}
\titleformat{\subsubsection}[runin]{\itshape}{\bfseries\upshape\thesubsubsection.}{1ex}{}[.---]

\titlespacing*{\section}{0pt}{*3}{*1}
\titlespacing*{\subsection}{0pt}{*3}{*1}
\titlespacing*{\subsubsection}{0pt}{*0.8}{*0}

\makeatletter
\newcommand{\neutralize}[1]{\expandafter\let\csname c@#1\endcsname\count@}
\makeatother

\newenvironment{theocopy}[1]
{\renewcommand{\theequation}{\ref{#1}}%
	\neutralize{equation}\phantomsection
	\begin{theo}}
	{\end{theo}}

\newenvironment{lemmacopy}[1]
{\renewcommand{\theequation}{\ref{#1}}%
	\neutralize{equation}\phantomsection
	\begin{lemma}}
	{\end{lemma}}

\renewbibmacro{in:}{}

\renewbibmacro*{volume+number+eid}{%
	\printfield{volume}%
	\setunit*{\addnbspace}
	\printfield{number}%
	\setunit{\addcomma\space}%
	\printfield{eid}}

\DeclareFieldFormat[article]{volume}{\textbf{#1}\space}
\DeclareFieldFormat[article]{number}{\mkbibparens{#1}}

\DeclareFieldFormat{journaltitle}{#1,}
\DeclareFieldFormat[thesis]{title}{\mkbibemph{#1}\addperiod}
\DeclareFieldFormat[article, unpublished, thesis]{title}{\mkbibemph{#1},}
\DeclareFieldFormat[book]{title}{\mkbibemph{#1}\addperiod}
\DeclareFieldFormat[unpublished]{howpublished}{#1, }

\DeclareFieldFormat{pages}{#1}

\DeclareFieldFormat[article]{series}{Ser.~#1\addcomma}

\setlist{topsep=3pt,itemsep=3pt}

\makeatletter
\newenvironment{qedequation*}{%
	\pushQED{\qed}%
	\mathdisplay@push
	\st@rredtrue \global\@eqnswfalse
	\mathdisplay{equation*}%
}{%
	\endmathdisplay{equation*}%
	\mathdisplay@pop
	\ignorespacesafterend
	\popQED\@endpefalse
}
\makeatother

\begin{document}
	\pagestyle{plain}
	
	\vspace*{-10pt}
	
	\maketitle
	
	\begin{abstract}
		The probabilistic method is a technique for proving combinatorial existence results by means of showing that a randomly chosen object has the desired properties with positive probability. A particularly powerful probabilistic tool is the Lov\'asz Local Lemma \ep{the LLL for short}, which was introduced by Erd\H os and Lov\'asz in the mid-1970s. Here we develop a version of the LLL that can be used to prove the existence of continuous colorings. 
		We then give several applications in Borel and topological dynamics. 
		\begin{itemize}[wide]
			\item Seward and Tucker-Drob showed that every free Borel action $\G \acts X$ of a countable group $\G$ admits an equivariant Borel map $\pi \colon X \to Y$ to a free subshift $Y \subset 2^\G$. We give a new simple proof of this result.
			
			\item We show that for a countable group $\G$, $\Free(2^\G)$ is weakly contained, in the sense of Elek, in every free continuous action of $\G$ on a zero-dimensional Polish space. This fact is analogous to the theorem of Ab\'ert and Weiss for probability measure-preserving actions and has a number of consequences in continuous combinatorics. In particular, we deduce that a coloring problem admits a continuous solution on $\Free(2^\G)$ if and only if it can be solved on finite subgraphs of the Cayley graph of $\G$ by an efficient deterministic distributed algorithm \ep{this fact	was also proved independently and using different methods by Seward}. This establishes a formal correspondence between questions that have been studied independently in continuous combinatorics and in distributed computing.
		\end{itemize}
	\end{abstract}

	\section{Introduction}
	
	\subsection{A continuous version of the Lov\'asz Local Lemma}
	
	\subsubsection{Constraint satisfaction problems and the LLL}
	
	Suppose we wish to prove that an object with certain combinatorial properties exists. A possible way to achieve this is by showing that an object chosen \emph{at random} from some class has the desired properties with positive probability. This approach was pioneered by Erd\H os in the 1940s and has since become indispensable throughout combinatorics; 
	see the book \cite{AS} by Alon and Spencer for an introduction. An important probabilistic tool is the so-called \emphd{Lov\'asz Local Lemma} \ep{the \emphd{LLL} for short}. The LLL is particularly useful for proving the existence of colorings satisfying a given set of ``local'' constraints. Formally, we define \emph{constraint satisfaction problems} as follows:
	
	\begin{defn}
		Let $X$ be a set and let $k \in \N^+$. We identify $k$ with the $k$-element set $\set{0, 1, \ldots, k-1}$.
		
		\begin{itemize}[wide]
			\item A \emphd{$k$-coloring} of a set $S$ is a function $f \colon S \to k$.
			
			\item For a finite set $D \subseteq X$, an \emphd{$(X,k)$-constraint} \ep{or simply a \emphd{constraint} if $X$ and $k$ are clear from the context} with \emphd{domain} $D$ is a set $B \subseteq k^D$ of $k$-colorings of $D$. 
			We write $\dom(B) \defeq D$.
			
			\item A $k$-coloring $f \colon X \to k$ \emphd{violates} a constraint $B$ with domain $D$ if the restriction of $f$ to $D$ is in $B$, and \emphd{satisfies} $B$ otherwise.
			
			\item A \emphd{constraint satisfaction problem} \ep{a \emphd{CSP} for short} $\B$ on $X$ with range $k$, in symbols $\B \colon X \to^? k$, is a set of $(X,k)$-constraints.
			
			\item A \emphd{solution} to a CSP $\B \colon X \to^? k$ is a $k$-coloring $f \colon X \to k$ that satisfies every constraint $B \in \B$.
		\end{itemize}
	\end{defn}

	In other words, each constraint $B \in \B$ in a CSP $\B \colon X \to^? k$ is interpreted as a set of finite ``forbidden patterns'' that are not allowed to appear in a solution $f \colon X \to k$. The LLL provides a simple probabilistic condition that guarantees that a given CSP has a solution. Fix a CSP $\B \colon X \to^? k$. For each $B \in \B$, the \emphd{probability} of $B$ is the quantity $\P[B]$ defined by
	\[
		\P[B] \,\defeq\, \frac{|B|}{k^{|\dom(B)|}} \,=\, \text{the probability that $B$ is violated by uniformly random $f \colon X \to k$}.
	\]
	The \emphd{neighborhood} of $B$ is the set
	\[
		\Nbhd(B) \,\defeq\, \set{B' \in \B \,:\, B'\neq B \text{ and } \dom(B') \cap \dom(B) \neq \0}.
	\]
	The LLL invokes the parameters $\pr(\B) \defeq \sup_{B \in \B} \P[B]$ and $\de(\B) \defeq \sup_{B \in \B} |\Nbhd(B)|$. 
	
	\begin{theo}[{{Lov\'asz Local Lemma}; Erd\H os--Lov\'asz \cite{EL}}]\label{theo:LLL}
		If $\B$ is a CSP such that
		\begin{equation}\label{eq:LLL}
			e \cdot \pr(\B) \cdot (\de(\B) + 1) \,\leq\, 1,
		\end{equation}
		where $e = 2.71\ldots$ is the base of the natural logarithm, then $\B$ has a solution.
	\end{theo}

	The LLL is often stated in the case when $\B$ is finite. However, a straightforward compactness argument shows that Theorem~\ref{theo:LLL} holds for infinite $\B$ as well \ep{see, e.g., \cite[proof of Theorem 5.2.2]{AS}}.
	
	\subsubsection{Continuous colorings}\label{subsubsec:cont_LLL}
	
	In this paper we are interested in the following question:
	
	\begin{ques}[{Continuous LLL}]\label{ques:cont_ques}
		Suppose $X$ is a zero-dimensional Polish space. What LLL-style conditions guarantee that a CSP $\B \colon X \to^? k$ has a \emph{continuous} solution $f \colon X \to k$? 
	\end{ques}
	
	Recall that a topological space is \emphd{Polish} if it is separable and completely metrizable, and \emphd{zero-dimensional} if it has a base consisting of clopen sets. The restriction to zero-dimensional spaces $X$ in Question~\ref{ques:cont_ques} is natural since a continuous map $X \to k$ can be thought of as a partition of $X$ into clopen sets indexed by $0$, $1$, \ldots, $k-1$, so if we hope to find a continuous solution to $\B$, it is reasonable to assume that $X$ has ``many'' clopen subsets. 
	Questions in the spirit of Question~\ref{ques:cont_ques} have recently attracted attention due to their applications in dynamical systems and descriptive set theory. For a sample of related results, see \cite{MLLL, CGMPT, BerDist}. These questions form a part of the general area called \emphd{descriptive combinatorics}, which investigates combinatorial problems under a variety of topological or measure-theoretic regularity requirements. For more background, see the surveys \cite{KechrisMarks} by Kechris and Marks and \cite{Pikh_survey} by Pikhurko.

	In the context of Question~\ref{ques:cont_ques}, it is necessary to assume that the CSP $\B$ itself ``respects'' the topology on $X$ in an appropriate sense. To this end, we define \emph{continuous CSPs} as follows:
	
	\begin{defn}\label{defn:cont}
		Let $X$ be a zero-dimensional Polish space. A CSP $\B \colon X \to^? k$ is \emphd{continuous} if for every set $B$ of functions $\set{1,\ldots, n} \to k$ and for all clopen subsets $U_2$, \ldots, $U_n \subseteq X$, the following set is clopen:
		\[
		\set{x_1 \in X \,:\, \text{there are $x_2 \in U_2$, \ldots, $x_n \in U_n$ such that $x_1$, \ldots, $x_n$ are distinct and $B(x_1, \ldots, x_n) \in \B$}}.
		\]
		Here $B(x_1, \ldots, x_n) \defeq \set{\phi \circ \iota \,:\, \phi \in B}$, where $\iota \colon \set{x_1, \ldots, x_n} \to \set{1, \ldots, n}$ is given by $x_i \mapsto i$.
	\end{defn}

	
	
	Conley, Jackson, Marks, Seward, and Tucker-Drob \cite[Theorem~1.6]{CJMST-D} constructed examples showing that the standard LLL condition \eqref{eq:LLL} is not sufficient to guarantee the existence of a {Borel}---let alone continuous---solution. In contrast to this, we prove that a certain strengthening of \eqref{eq:LLL} does yield continuous solutions. In addition to $\pr(\B)$ and $\de(\B)$, we consider 
	two more parameters associated to a CSP $\B \colon X \to^? k$. Namely, we define the \emphd{maximum vertex-degree} $\vdeg(\B)$ of $\B$ as
	\[
		\vdeg(\B) \,\defeq\, \sup_{x \in X} |\set{B\in \B \,:\, x \in \dom(B)}|,
	\]
	and let the \emphd{order} $\ord(\B)$ of $\B$ be $\ord(\B) \defeq \sup_{B \in \B} |\dom(B)|$. Note that $\de(\B) \leq (\vdeg(\B) - 1) \ord(\B)$. 
	
	\begin{theo}\label{theo:cont_LLL}
		Let $\B \colon X \to^? k$ be a continuous CSP on a zero-dimensional Polish space $X$. If
		\begin{equation}\label{eq:bound}
			\pr(\B) \cdot \vdeg(\B)^{\ord(\B)} \,<\, 1,
		\end{equation}
		then $\B$ has a continuous solution $f \colon X \to k$.
	\end{theo}
	
	Note that in the setting of Theorem~\ref{theo:cont_LLL}, if $\pr(\B) > 0$, then in fact $\pr(\B) \geq k^{-\ord(\B)}$. Thus, Theorem~\ref{theo:cont_LLL} is only useful if $k$ is relatively large \ep{namely $k > \vdeg(\B)$}.
	
	
	We prove Theorem~\ref{theo:cont_LLL} in \S\ref{sec:proof_LLL} using the \emph{method of conditional probabilities}---a standard derandomization technique in computer science. This connection to computer science is not coincidental: results and methods in descriptive combinatorics often mirror those in \emphd{distributed computing}, i.e., the area concerned with problems that can be solved efficiently by a decentralized network of processors. For example, an argument similar to our proof of Theorem~\ref{theo:cont_LLL} was involved in Fischer and Ghaffari's breakthrough work on distributed algorithms for the LLL \cite[Theorem 3.5]{FG}.
	
	Another relevant result in distributed computing is due to Brandt, Grunau, and Rozho\v{n} \cite{expLLL}, who recently developed an efficient deterministic distributed algorithm for finding solutions to CSPs under the condition $\pr 2^{\de} < 1$ \ep{in the special case when $\vdeg \leq 3$, such an algorithm was devised earlier by Brandt, Maus, and Uitto \cite{BMU}}. 
	In \cite{BerDist}, the author established a series of general results that allow using efficient distributed algorithms to obtain colorings with desirable regularity properties \ep{such as continuity, measurability, etc.}. In particular, \cite[Theorem~2.13]{BerDist} implies that under suitable assumptions the condition $\pr 2^{\de} < 1$ is also sufficient to produce continuous solutions. While in general neither of the bounds $\pr 2^{\de} < 1$ and $\pr \cdot \vdeg^\ord < 1$ \ep{that is, \eqref{eq:bound}} implies the other, in practice one often estimates $\de$ using the inequality $\de \leq (\vdeg-1) \ord$, which makes the latter bound more widely applicable, especially when $\ord$ is much smaller than $\vdeg$.
	
	Remarkably, the bound $\pr 2^{\de} < 1$ in the Brandt--Grunau--Rozho\v{n} result is sharp: there is no such efficient distributed algorithm that finds solutions to CSPs if the bound is relaxed to $\pr 2^{\de} \leq 1$. This follows from the analysis of the so-called \emph{sinkless orientation problem} performed in the randomized setting by Brandt et al. \cite{lowerbound} and extended to the deterministic setting by Chang, Kopelowitz, and Pettie \cite{CKP}. This sharpness result has a counterpart in descriptive combinatorics. Namely, suppose $\de \in \N$ and let $G$ be a \emphd{$\de$-regular} graph, meaning that every vertex of $G$ is incident to exactly $\de$ edges. An orientation of $G$ is \emphd{sinkless} if the outdegree of every vertex is at least $1$. A sinkless orientation of $G$ can be naturally encoded as a solution to a CSP $\B_{\text{sinkless}} = \set{B_x}_{x \in V(G)} \colon E(G) \to^? 2$. Here the color of each edge $e \in E(G)$ indicates the direction in which $e$ is oriented, and $B_x$ for $x \in V(G)$ is the constraint with domain $\dom(B_x) = \set{e \in E(G) \,:\, \text{$e$ is incident to $x$}}$ that requires the outdegree of $x$ to be at least $1$. It is easy to see that $\de(\B_{\text{sinkless}}) = \de$ and $\pr(\B_{\text{sinkless}}) = 1/2^\de$. However, Thornton \cite[Theorem 3.5]{Thor} used the determinacy method of Marks \cite{Marks} to construct, for any given $\de \in \N$, a Borel $\de$-regular graph $G$ that does not admit a Borel sinkless orientation. Note that since $\vdeg(\B_{\text{sinkless}}) = 2$ and $\ord(\B_{\text{sinkless}}) = \de$, this also serves as a sharpness example for Theorem~\ref{theo:cont_LLL}.
	
	We shall resume the discussion of distributed algorithms in \S\ref{subsubsec:LOCAL}, where we describe one of the consequences derived using Theorem~\ref{theo:cont_LLL}, namely that for certain types of coloring problems, a continuous solution exists \emph{if and only if} the problem can be solved by an efficient distributed algorithm.
	
	\subsubsection{Borel colorings}
	
	Sometimes we only wish to find a Borel solution instead of a continuous one. Recall that a \emphd{standard Borel space} is a set $X$ equipped with a $\sigma$-algebra $\mathfrak{B}(X)$ of \emphd{Borel sets} generated by a Polish topology on $X$. We say that a Polish topology on a standard Borel space $X$ is \emphd{compatible} if it generates $\mathfrak{B}(X)$. If $X$ is a standard Borel space, then the set $\fins{X}$ of all finite subsets of $X$ also carries a natural standard Borel structure. Since every $(X,k)$-constraint can be viewed as a finite subset of $\fins{X \times k}$, we may speak of \emphd{Borel CSPs} $\B \colon X \to^? k$, i.e., Borel sets $\B \subseteq \fins{\fins{X \times k}}$ of $(X,k)$-constraints. The following is an immediate corollary of Theorem~\ref{theo:cont_LLL}:

	\begin{corl}
		Let $\B \colon X \to^? k$ be a Borel CSP on a standard Borel space $X$. If
		\[
		\pr(\B) \cdot \vdeg(\B)^{\ord(\B)} \,<\, 1,
		\]
		then $\B$ has a Borel solution $f \colon X \to k$.
	\end{corl}
	\begin{scproof}
		For a set $B$ of functions $\set{1,\ldots, n} \to k$, write $x_1 \sim_B (x_2, \ldots, x_n)$ if $x_1$, \ldots, $x_n$ are distinct and $B(x_1, \ldots, x_n) \in \B$. Since $\vdeg(\B) < \infty$, the Luzin--Novikov theorem \cite[Theorem 18.10]{KechrisDST} yields a finite sequence of partial Borel maps $h_{B,i} \colon X \rightharpoonup X^{n-1}$, $1 \leq i \leq \vdeg(\B)(n-1)!$, such that $x_1 \sim_B (x_2, \ldots, x_n)$ if and only if $(x_2, \ldots, x_n) = h_{B, i}(x_1)$ for some $i$. \ep{The $(n-1)!$ factor arises from the fact that a constraint with domain $\set{x_1, \ldots, x_n}$ could, in principle, force $x_1$ to be $\sim_B$-related to every permutation of $\set{x_2, \ldots, x_n}$.} It follows from standard results in descriptive set theory \cite[\S13]{KechrisDST} that there is a compatible zero-dimensional Polish topology $\tau$ on $X$ with respect to which all the maps $h_{B,i}$ are continuous and defined on clopen sets. Then $\B$ is continuous with respect to $\tau$, so, by Theorem~\ref{theo:cont_LLL}, $\B$ has a  $\tau$-continuous \ep{hence Borel} solution. Alternatively, it is straightforward to check directly that the proof of Theorem~\ref{theo:cont_LLL} given in \S\ref{sec:proof_LLL} goes through in the Borel setting with the words ``continuous'' and ``clopen'' replaced everywhere by ``Borel.''
	\end{scproof}

	\subsection{Applications in dynamics} 
	
	
	\subsubsection{A simple proof of the Seward--Tucker-Drob theorem}
	
	Throughout the rest of this paper, $\G$ denotes a countably infinite discrete group with identity element~$\mathbf{1}$. By an ``action'' of $\G$ we always mean a left action. Our first application of Theorem~\ref{theo:cont_LLL} is a simple proof of the following result of Seward and Tucker-Drob:
	
	\begin{theo}[{Seward--Tucker-Drob \cite{STD}}]\label{theo:STD}
		If $\G \acts X$ is a free Borel action of $\G$ on a standard Borel space $X$, then there is a $\G$-equivariant Borel map $\pi \colon X \to Y$, where $Y \subset 2^\G$ is a free subshift.
	\end{theo}
	
	Let us recall the terminology used in the statement of Theorem~\ref{theo:STD}. A~\emphd{$\G$-space} is a topological space $X$ equipped with a continuous action $\G \acts X$. 
	The product space $k^\G$ of all $k$-colorings $\G \to k$ of $\G$ is a compact zero-dimensional Polish space, and it becomes a $\G$-space under the action $\G\acts k^\G$ given by\footnote{Alternatively, we could define $(\gamma \cdot x)(\delta)$ to be $x(\gamma^{-1}\delta)$. The two definitions yield isomorphic structures, but we prefer to avoid the use of inverses.} 
	\[(\gamma \cdot x)(\delta) \,\defeq\, x(\delta \gamma) \quad \text{for all } x \in k^\G \text{ and } \gamma,\, \delta \in \G.\] The $\G$-spaces of the form $k^\G$ are called \emphd{Bernoulli shifts}, or simply \emphd{shifts}. The \emphd{free part} of a $\G$-space $X$ is the set $\Free(X) \defeq \set{x \in X \,:\, \Stab_\G(x) = \set{\mathbf{1}}}$ equipped with the subspace topology and the induced action of $\G$ \ep{here $\Stab_\G(x)$ denotes the stabilizer of $x$}. In other words, $\Free(X)$ is the largest $\G$-invariant subspace of $X$ on which $\G$ acts freely. If $X$ is a Polish $\G$-space, then $\Free(X)$ is a $G_\delta$ subset of $X$, and hence it is also Polish \cite[Theorem 3.11]{KechrisDST}. A closed $\G$-invariant subset of $k^\G$ is called a \emphd{subshift}, and we say that a subshift $X \subseteq k^\G$ is \emphd{free} if $X \subseteq \Free(k^\G)$, i.e., if $\G$ acts freely on $X$.

	Even the \emph{existence} of a nonempty free subshift for an arbitrary countable group $\G$ is far from obvious. It was established by Gao, Jackson, and Seward \cite{GJS1} with a rather technical construction that was further explored in \cite{GJS2}. The proof of Theorem~\ref{theo:STD} due to Seward and Tucker-Drob develops the ideas of Gao, Jackson, and Seward further and is similarly quite involved. However, it turns out that questions about free subshifts are well-suited for LLL-based approaches. This fact was first observed by Aubrun, Barbieri, and Thomass\'e \cite{ABT}, who gave a short and simple LLL-based alternative construction of a nonempty free subshift $X \subset 2^\G$. Roughly speaking, their method was to define $X$ as the set of all colorings $\G \to 2$ satisfying certain constraints, and then to show that $X \neq \0$ using the LLL. This technique is quite flexible and allows constructing free subshifts with various additional properties. For example, in \cite{LSS} the LLL is applied to construct free subshifts that are not just nonempty but ``large'' {in terms of Hausdorff dimension and entropy}.
	
	Unfortunately, the approach of \cite{ABT, LSS} cannot prove Theorem~\ref{theo:STD}, since it invokes a version of the LLL that does not generally yield Borel solutions. 
	\ep{Actually, \cite{ABT, LSS} rely on the so-called \emph{General LLL} \cite[Lemma 5.1.1]{AS}, which is a strengthening of Theorem~\ref{theo:LLL} that, in general, does not even yield \emph{measurable} solutions \cite[Theorem 7.1]{MLLL}.} In \S\ref{sec:main_proof}, we show that Theorem~\ref{theo:STD} can nevertheless be established with a simple probabilistic argument---namely with the help of Theorem~\ref{theo:cont_LLL}.
	
	\subsubsection{Topological Ab\'ert--Weiss theorem}
	
	Inspired by the analogous notions for measure-preserving actions \ep{which were in turn modeled after similar concepts in representation theory}, Elek \cite{Elek} introduced the relations of \emph{weak containment} and \emph{weak equivalence} on the class of zero-dimensional Polish $\G$-spaces. \ep{Technically, Elek only considered \emph{compact} zero-dimensional $\G$-spaces, but the same definitions can be applied verbatim to non-compact spaces as well.} Let $k \geq 1$ be an integer. A~\emphd{$k$-pattern} is a partial map $p \colon \G \rightharpoonup k$ whose domain is a finite subset of $\G$. Given an action $\G \acts X$ and a $k$-coloring $f \colon X \to k$, we say that a $k$-pattern $p$ \emphd{occurs} in $f$ if there is a point $x \in X$ such that $f(\gamma \cdot x) = p(\gamma)$ for all $\gamma \in \dom(p)$. 
	For a finite subset $F \subset \G$, we let $\Pat_F(X, f)$ denote the set of all $k$-patterns $p \colon F \to k$ with domain $F$ that occur in $f$. 
	
	\begin{defn}\label{defn:weak}
		Let $X$ and $Y$ be zero-dimensional Polish $\G$-spaces. We say that $X$ is \emphd{weakly contained} in $Y$, in symbols $X \preccurlyeq Y$, if given any $k \in \N^+$, a finite subset $F \subset \G$, and a continuous $k$-coloring $f \colon X \to k$, there is a continuous $k$-coloring $g \colon Y \to k$ such that $\Pat_F(Y, g) = \Pat_F(X, f)$. If $X \preccurlyeq Y$ and $Y \preccurlyeq X$, then we say that $X$ and $Y$ are \emphd{weakly equivalent} and write $X \simeq Y$.
	\end{defn}

	As mentioned earlier, Definition~\ref{defn:weak} was introduced \ep{for compact $\G$-spaces} by Elek in \cite{Elek}. For minimal actions of the group $\Z$, weak equivalence \ep{under the name of \emph{weak approximate conjugacy}} was considered previously by Lin and Matui in \cite{LM}.
	
	Among several other results, Elek proved that the pre-order of weak containment has a minimum element in the class of all nonempty free zero-dimensional Polish $\G$-spaces \cite[Theorem 2]{Elek}. In other words, Elek showed that there exists a free \ep{compact} zero-dimensional Polish $\G$-space $M$ such that $M \preccurlyeq X$ for every nonempty free zero-dimensional Polish $\G$-space $X$ \ep{it is easy to check that Elek's argument does not need $X$ to be compact}. We show that, except for the compactness requirement, one can actually take $M$ to be the free part of the Bernoulli shift $2^\G$:
	
	\begin{theo}\label{theo:top_AW}
		If $X$ is a nonempty free zero-dimensional Polish $\G$-space, then $\Free(2^\G) \preccurlyeq X$.
	\end{theo}

	Theorem~\ref{theo:top_AW} is a topological counterpart to the ergodic-theoretic result of Ab\'ert and Weiss \cite{AW}, namely that the Bernoulli shift $2^\G$ is weakly contained \ep{in the sense of Kechris \cite{K_book}} in each almost everywhere free probability measure-preserving action of $\G$. The proof of Theorem~\ref{theo:top_AW} is given in \S\ref{subsec:proof_top_AW}. It is an elaboration of our proof of Theorem~\ref{theo:STD}, leveraging the fact that Theorem~\ref{theo:cont_LLL} yields continuous \ep{and not just Borel} solutions.

	\subsection{Consequences in continuous combinatorics}
	
	
	The main motivation for this work comes from the area of \emphd{continuous combinatorics}, which studies the behavior of combinatorial notions---such as graph colorings, matchings, etc.---under additional continuity constraints. For example, suppose that $G$ is a graph whose vertex set $V(G)$ is a zero-dimensional Polish space. A typical problem in continuous combinatorics is to determine the \emphd{continuous chromatic number} $\chi_c(G)$ of $G$, i.e., the least $k$ for which there exists a continuous $k$-coloring $f \colon V(G) \to k$ satisfying $f(x) \neq f(y)$ whenever vertices $x$ and $y$ are adjacent \ep{such colorings are called \emphd{proper}}. 
	
	In \cite{Abelian}, Gao, Jackson, Krohne, and Seward initiated the systematic study of continuous combinatorics of countable group actions and performed a detailed analysis in the case $\G = \Z^d$. In particular, they completely characterized combinatorial problems that can be solved continuously on the space $\Free(2^{\G})$ for $\G \in \set{\Z, \Z^2}$ by reducing them to certain questions about finite graphs. Here we continue this line of research and extend it to the case of $\G$-spaces for arbitrary countably infinite groups $\G$.
	
	Some of our results in this section, specifically Theorems~\ref{theo:universal} and \ref{theo:LOCAL}, 
	were obtained independently by Seward using the techniques from \cite{GJS1, GJS2, STD} \ep{personal communication}. 
	
	\subsubsection{Universality of the shift}
	
	We say that a coloring $f \colon X \to k$ is \emphd{$\Pat$-avoiding}, where $X$ is a $\G$-space and $\Pat$ is a set of $k$-patterns, if no pattern $p \in \Pat$ occurs in $f$. As a side remark, we note that continuous $\Pat$-avoiding colorings of $\G$-spaces have a natural meaning from the standpoint of topological dynamics. Specifically, viewing $\G$ itself as a discrete $\G$-space under the left multiplication action $\G \acts \G$, we can consider the set $\mathsf{Av}(\Pat) \subseteq k^\G$ of all $\Pat$-avoiding $k$-colorings of $\G$, for a given finite set $\Pat$ of $k$-patterns. The set $\mathsf{Av}(\Pat)$ is closed and $\G$-invariant, and it is called a \emphd{subshift of finite type} \ep{``finite'' because $\Pat$ is finite}. If $X$ is a $\G$-space, then there is a natural one-to-one correspondence
	\[
	\set{\text{$\Pat$-avoiding continuous colorings $X \to k$}} \quad \longleftrightarrow \quad \set{\text{$\G$-equivariant continuous maps $X \to \mathsf{Av}(\Pat)$}},
	\]
	where each 
	$\Pat$-avoiding continuous coloring $f \colon X \to k$ gives rise to the so-called \emphd{coding map} $\pi_f \colon X \to \mathsf{Av}(\Pat)$ given by $\pi_f(x)(\gamma) \defeq f(\gamma \cdot x)$ for all $x \in X$ and $\gamma \in \G$.
	In view of this correspondence, studying continuous colorings that avoid finite sets of patterns is equivalent to studying equivariant continuous maps to subshifts of finite type. The following is an immediate consequence of Theorem~\ref{theo:top_AW}:
	
	\begin{theo}\label{theo:universal}
		Let $\Pat$ be a finite set of $k$-patterns. The following statements are equivalent.
		\begin{enumerate}[label=\ep{\normalfont\arabic*}]
			\item\label{item:shift} There is a continuous $\Pat$-avoiding $k$-coloring of $\Free(2^\G)$.
			
			\item\label{item:all} Every free zero-dimensional Polish $\G$-space admits a continuous $\Pat$-avoiding $k$-coloring. 
		\end{enumerate}
	\end{theo}
	\begin{scproof}
		Implication \ref{item:all} $\Longrightarrow$ \ref{item:shift} is obvious, while \ref{item:shift} $\Longrightarrow$ \ref{item:all} is given by Theorem~\ref{theo:top_AW}.
	\end{scproof}

	Informally, Theorem~\ref{theo:universal} says that of all the free zero-dimensional Polish $\G$-spaces, it is the hardest to solve combinatorial problems continuously on $\Free(2^\G)$. Here is just one specific instance of this phenomenon. Let $S \subset \G$ be finite. 
	The \emphd{Schreier graph} of a $\G$-space $X$ corresponding to $S$ is the \ep{simple undirected} graph $G(X, S)$ with vertex set $X$ where two distinct vertices $x$, $y$ are adjacent if and only if $y = \sigma \cdot x$ for some $\sigma \in S \cup S^{-1}$. A consequence of Theorem~\ref{theo:universal} is that the Schreier graph of $\Free(2^\G)$ has the largest continuous chromatic number among all Schreier graphs of free zero-dimensional Polish $\G$-spaces:

	\begin{corl}\label{corl:color}
		Let $S$ be a finite subset of $\G$. If $X$ is a free zero-dimensional Polish $\G$-space, then \[\chi_c(G(X, S)) \,\leq\, \chi_c(G(\Free(2^\G), S)).\] 
	\end{corl}
	\begin{scproof}
		Set $k \defeq \chi_c(G(\Free(2^\G), S))$ and apply Theorem~\ref{theo:universal} with $\Pat \defeq \set{p_{i, \sigma} \,:\, 0 \leq i < k, \sigma \in S \setminus \set{\mathbf{1}}}$, where for each $i$ and $\sigma$, $p_{i, \sigma}$ is the $k$-pattern with domain $\set{\mathbf{1}, \sigma}$ that sends both $\mathbf{1}$ and $\sigma$ to $i$.
	\end{scproof}
	
	\subsubsection{Reduction to finite graphs}
	
	In our remaining results, we reduce problems about continuous colorings to questions about colorings of finite graphs. To state them, we require a few definitions. Let $S \subset \G$ be a finite set. An \emphd{$S$-labeled graph} is a simple undirected graph $G$ equipped with a \emphd{labeling map} $\lambda$ that assigns to each \ep{ordered} pair $(x,y)$ of adjacent vertices a group element $\lambda(x,y) \in S \cup S^{-1}$ so that $\lambda(y,x) = \lambda(x,y)^{-1}$. Note that a vertex $x$ may have multiple neighbors $y$ with the same $\lambda(x,y)$. For a subset $U \subseteq V(G)$, we let $G[U]$ denote the \emphd{subgraph} of $G$ \emphd{induced} by $U$, i.e., the $S$-labeled graph with vertex set $U$ whose adjacency relation and labeling map are inherited from $G$.
	
	Schreier graphs of free $\G$-spaces are natural examples of $S$-labeled graphs, with $\lambda(x,y)$ being the unique element $\sigma \in S \cup S^{-1}$ such that $y = \sigma \cdot x$. When $\G$ itself is viewed as a discrete $\G$-space under the left multiplication action $\G \acts \G$, the $S$-labeled Schreier graph $G(\G, S)$ is called the \emphd{Cayley graph} of $\G$ corresponding to $S$. Note that the graph $G(\G, S)$ is connected if and only if $S$ generates $\G$. 
	For a subset $F \subseteq \G$, we use $G(F, S) \defeq G(\G,S)[F]$ to denote the subgraph of $G(\G, S)$ induced by $F$.
	
	A \emphd{homomorphism} from an $S$-labeled graph $G$ to an $S$-labeled graph $H$ is a map $\phi \colon V(G) \to V(H)$ such that if $x$, $y \in V(G)$ are adjacent in $G$, then $\phi(x)$, $\phi(y)$ are adjacent in $H$ and $\lambda(\phi(x), \phi(y)) = \lambda(x, y)$. Let $F \subset \G$ be a finite set and let $p \colon F \to k$ be a $k$-pattern. We say that $p$ is \emphd{$S$-connected} if the graph $G(F,S)$ is connected. Given an $S$-labeled graph $G$ and a coloring $f \colon V(G) \to k$, we say that an $S$-connected $k$-pattern $p \colon F \to k$ \emphd{occurs} in $f$ if there is a homomorphism $\phi \colon F \to V(G)$ from $G(F, S)$ to $G$ such that $f \circ \phi = p$. When $G$ is the Schreier graph $G(X, S)$ of a free $\G$-space $X$, this notion coincides with our previous definition, since the only homomorphisms from $G(F, S)$ to $G(X, S)$ are the ones of the form $F  \to X \colon \gamma \mapsto \gamma \cdot x$ for some $x \in X$ \ep{here we use that $p$ is $S$-connected}. Given a finite set $\Pat$ of $S$-connected $k$-patterns, we say that a coloring $f \colon V(G) \to k$ of an $S$-labeled graph $G$ is \emphd{$\Pat$-avoiding} if none of the patterns in $\Pat$ occur in $f$.
	
	Consider the standard generating set $S \defeq \set{(1,0), (0,1)}$ for the group $\Z^2$. In \cite[Theorem~5.5]{Abelian}, Gao, Jackson, Krohne, and Seward constructed an explicit countable family $\mathscr{H}$ of finite $S$-labeled graphs such that the following statements are equivalent for any finite set $\Pat$ of $S$-connected $k$-patterns:
	\begin{itemize}
		\item $\Free(2^{\Z^2})$ admits a continuous $\Pat$-avoiding $k$-coloring;
		
		\item there is a graph in $\mathscr{H}$ that admits a $\Pat$-avoiding $k$-coloring;
		
		\item all but finitely many graphs in $\mathscr{H}$ admit $\Pat$-avoiding $k$-colorings.
	\end{itemize}
	In other words, to determine whether the \ep{infinite} $\Z^2$-space $\Free(2^{\Z^2})$ has a continuous $\Pat$-avoiding $k$-coloring, one simply needs to check if the \ep{finite} graphs in $\mathscr{H}$ admit $\Pat$-avoiding $k$-colorings. This can be seen as a ``compactness theorem'' for continuous colorings of $\Free(2^{\Z^2})$. Gao, Jackson, Krohne, and Seward call \cite[Theorem~5.5]{Abelian} the ``Twelve Tiles Theorem,'' since each graph in $\mathscr{H}$ is obtained from twelve pieces---``tiles''---glued to each other according to certain rules.
	
	We obtain 
	an analogous result for arbitrary countable groups $\G$:
	
	\begin{theo}\label{theo:compact}
		In the setting of Theorem~\ref{theo:universal}, assume that $S \subset \G$ is a finite set such that the $k$-patterns in $\Pat$ are $S$-connected. There is an explicit countable family $\mathscr{H}$ of finite $S$-labeled graphs \ep{see \S\ref{subsec:tiles} for the definition} such that statements \ref{item:shift} and \ref{item:all} are also equivalent to:
		\begin{enumerate}[label=\ep{\normalfont\arabic*}]\setcounter{enumi}{2}
			\item\label{item:compact1} There is a graph in $\mathscr{H}$ that admits a $\Pat$-avoiding $k$-coloring.
			
			\item\label{item:compact_many} All but finitely many graphs in $\mathscr{H}$ admit $\Pat$-avoiding $k$-colorings.
		\end{enumerate}
	\end{theo}
	
	We construct the family $\mathscr{H}$ and prove Theorem~\ref{theo:compact} in \S\ref{subsec:tiles}.
	
	\subsubsection{\LOCAL algorithms}\label{subsubsec:LOCAL}
	
	Our final result establishes a precise connection between continuous combinatorics and distributed computing:
	
	\begin{theo}\label{theo:LOCAL}
		In the setting of Theorem~\ref{theo:universal}, assume that $S \subset \G$ is a finite set such that the $k$-patterns in $\Pat$ are $S$-connected. Then statements \ref{item:shift}--\ref{item:compact_many} are also equivalent to:
		\begin{enumerate}[label=\ep{\normalfont\arabic*}]\setcounter{enumi}{4}
			\item\label{item:LOCAL} There is a deterministic distributed algorithm in the \LOCAL model that, given an $n$-vertex $S$-labeled subgraph $G$ of $G(\G, S)$, in $O(\log^\ast n)$ rounds outputs a $\Pat$-avoiding $k$-coloring of $G$.
		\end{enumerate}
	\end{theo}
	
	Here $\log^\ast n$ denotes the \emphd{iterated logarithm} of $n$, i.e., the number of times the logarithm function must be iteratively applied to $n$ before the result becomes at most $1$. 
	
	Statement \ref{item:LOCAL} in Theorem~\ref{theo:LOCAL} refers to the \LOCAL model of distributed computation, which was introduced by Linial in \cite{Linial}. 
	For a comprehensive introduction to this model, see the book \cite{BE} by Barenboim and Elkin. The \LOCAL model operates on an $n$-vertex graph $G$. Here we think of $G$ as representing a decentralized communication network where each vertex plays the role of a processor and edges represent communication links. The computation proceeds in \emph{rounds}. During each round, the vertices first perform some local computations and then synchronously broadcast messages to all their neighbors. After a number of rounds, every vertex must output a color, and the resulting coloring of $V(G)$ is considered to be the output of the algorithm. The efficiency of such an algorithm is measured by the number of communication rounds required.
	
	An important feature of the \LOCAL model is that every vertex of $G$ is executing the same algorithm. Therefore, to make this model nontrivial, the vertices must be given a way of breaking symmetry. In the \emph{deterministic} variant of the \LOCAL model, this is achieved by assigning a unique identifier $\mathrm{Id}(x) \in \set{1,\ldots, n}$ to every vertex $x \in V(G)$. 
	The identifier assigned to a vertex $x$ is treated as part of $x$'s input; that is, $x$ ``knows'' what its own identifier is initially and can communicate this information to its neighbors. When we say that a deterministic \LOCAL algorithm \emph{solves} a coloring problem $\Pat$ on a given class $\mathscr{G}$ of finite graphs, we mean that the coloring it outputs on any graph from $\mathscr{G}$ is a valid solution to $\Pat$, regardless of the way the identifiers are assigned. The word ``deterministic'' distinguishes this model from the \emph{randomized} version, where the vertices are allowed to generate sequences of random bits. In this paper we will only be concerned with deterministic algorithms.
	
	
	If $x$ and $y$ are two vertices whose graph distance in $G$ is $T$, then no information from $y$ can reach $x$ in fewer than $T$ communication rounds \ep{this explains the name ``\LOCAL''}. Conversely, in $T$ rounds every vertex can collect all the data present at the vertices at distance at most $T$ from it. Thus, a $T$-round \LOCAL algorithm may be construed simply as a function that, given the structure of the radius-$T$ ball around $x$ \ep{including the assignment of the identifiers to its vertices}, outputs $x$'s color \cite[\S{}4.1.2]{BE}. 
	
	In general, the input graph $G$ may possess some additional structure \ep{such as an orientation, a fixed coloring of the vertices, etc.}. For example, in Theorem~\ref{theo:LOCAL} we consider \LOCAL algorithms operating on \emph{$S$-labeled} graphs $G$. This means that the labels on the edges of $G$ form part of the problem's input, and each vertex can discover the labels of the edges in its radius-$T$ ball in $T$ communication rounds.
	
	The formal equivalence between general classes of problems in continuous combinatorics and in distributed computing given by Theorem~\ref{theo:LOCAL} explains the parallels between specific results in these two areas. For example, suppose $\G = \Z^2$ and let $S \defeq \set{(1,0), (0,1)}$. Among numerous other results, Gao, Jackson, Krohne, and Seward proved in \cite{Abelian} that the continuous chromatic number of $G(\Free(2^{\Z^2}), S)$ is $4$, 
	and also that there is no algorithm for deciding, given a finite set $\Pat$ of $k$-patterns, whether $\Free(2^{\Z^2})$ admits a continuous $\Pat$-avoiding $k$-coloring. 
	In \cite{grids}, Brandt et al.~established analogous results for distributed algorithms on $n\times n$ grid graphs: proper $k$-colorings of such graphs can be computed by an $O(\log^\ast n)$-round \LOCAL algorithm for $k \geq 4$ but not for $k = 3$, 
	and there is no decision procedure that determines, for a given finite set $\Pat$ of $k$-patterns, whether $\Pat$-avoiding $k$-colorings of such graphs can be found by a $O(\log^\ast n)$-round \LOCAL algorithm. 
	Theorem~\ref{theo:LOCAL} provides a general reason underlying this analogy.
	
	The connection between continuous combinatorics and distributed algorithms was observed recently by Elek \cite{Elek} and the author \cite{BerDist}.  
	In particular, implication \ref{item:LOCAL} $\Longrightarrow$ \ref{item:all} is a special case of \cite[Theorem~2.13]{BerDist}, which is a general result that provides a way to use efficient deterministic \LOCAL algorithms to obtain continuous colorings.
	Thus, we only have to prove \ref{item:all} $\Longrightarrow$ \ref{item:LOCAL} here, which is done in \S\ref{subsec:LOCAL} by utilizing the specific construction of the family of finite graphs $\mathscr{H}$ from Theorem~\ref{theo:compact}.
	
	\subsubsection*{{{Acknowledgments}}}
	
	I am indebted to Jan Greb\'ik, Stephen Jackson, V\'aclav Rozho\v{n}, Brandon Seward, and Zolt\'an Vidny\'{a}nszky for many insightful discussions and for helpful comments on the manuscript. I am also grateful to the anonymous referee for useful suggestions.
	
	\section{Preliminaries}\label{sec:prelim}
	
	We shall require a few basic facts about continuous graph combinatorics. These facts may be somewhat less well-known than their Borel counterparts, so we prove them here for completeness. \ep{The proofs are standard and essentially present in \cite[\S4]{KechrisSoleckiTodorcevic}.}
	
	Let $G$ be a graph. For a subset $S \subseteq V(G)$, $\Nbhd_G(S)$ denotes the \emphd{neighborhood} of $S$ in $G$, i.e., the set of all vertices that have a neighbor in $S$. For a vertex $x \in V(G)$, we write $\Nbhd_G(x) \defeq \Nbhd_G(\set{x})$. A graph $G$ is \emphd{locally finite} if $\Nbhd_G(x)$ is finite for every $x \in X$. A set $I \subseteq V(G)$ is \emphd{independent} if $I \cap \Nbhd_G(I) = \0$, i.e., if no two vertices in $I$ are adjacent. For a subset $U \subseteq V(G)$, we use $G[U]$ to denote the \emphd{subgraph} of $G$ \emphd{induced} by $U$, i.e., the graph with vertex set $U$ whose adjacency relation is inherited from $G$, and we write $G - U \defeq G [V(G) \setminus U]$. We say that $G$ is a \emphd{continuous graph} if $V(G)$ is a zero-dimensional Polish space and for every clopen set $U \subseteq V(G)$ its neighborhood $\Nbhd_G(U)$ is also clopen. \ep{This is analogous to Definition~\ref{defn:cont}.} Note that if $G$ is a continuous graph and $U \subseteq V(G)$ is a clopen set of vertices, then the subgraph $G[U]$ of $G$ induced by $U$ is also continuous. 
	
	\begin{lemma}\label{lemma:countable}
		Every locally finite continuous graph $G$ admits a partition $V(G) = \bigsqcup_{n = 0}^\infty I_n$ into countably many clopen independent sets.
	\end{lemma}
	\begin{scproof}
		Let $\set{U_n \,:\, n \in \N}$ be a countable base for the topology on $V(G)$ consisting of clopen sets. For each $n \in \N$, let $V_n \defeq U_n \setminus \Nbhd_G(U_n)$. By construction, each $V_n$ is independent and, since $G$ is continuous, clopen. Since $G$ is locally finite, each $x \in V(G)$ has an open neighborhood disjoint from $\Nbhd_G(x)$, and hence $\bigcup_{n = 0}^\infty V_n = V(G)$. It remains to make the sets disjoint by setting $I_n \defeq V_n \setminus (V_0 \cup \ldots \cup V_{n-1})$.
	\end{scproof}

	\begin{lemma}\label{lemma:max}
		Every locally finite continuous graph $G$ has a clopen maximal independent set $I \subseteq V(G)$.
	\end{lemma}
	\begin{scproof}
		Let $V(G) = \bigsqcup_{n = 0}^\infty I_n$ be a partition into countably many clopen independent sets given by Lemma~\ref{lemma:countable}. Define a sequence of clopen subsets $I_n' \subseteq I_n$ recursively by setting $I_0' \defeq I_0$ and $I_{n+1}' \defeq I_{n+1} \setminus \Nbhd_G(I_0' \sqcup \ldots \sqcup I_n')$ for all $n \in \N$. By construction, the set $I \defeq \bigsqcup_{n=0}^\infty I_n'$ is a maximal independent set in $G$. Since $G$ is continuous, the sets $I_n'$ are clopen, and hence $I$ is open. But the sets $I_n \setminus I_n'$ are also clopen, so $V(G) \setminus I = \bigsqcup_{n=0}^\infty (I_n \setminus I_n')$ is open as well, and hence $I$ is clopen, as desired. 
	\end{scproof}

	The \emphd{maximum degree} $\Delta(G)$ of a graph $G$ is defined by $\Delta(G) \defeq \sup_{x \in V(G)} |\Nbhd_G(x)|$.

	\begin{lemma}\label{lemma:coloring}
		If $G$ is a continuous graph of finite maximum degree $\Delta$, then $\chi_c(G) \leq \Delta + 1$.
	\end{lemma}
	\begin{scproof}
		We need to find a partition of $V(G)$ into $\Delta + 1$ clopen independent sets. To this end, we iteratively apply Lemma~\ref{lemma:max} to obtain a sequence $I_0$, \ldots, $I_\Delta$ where each $I_n$ is a clopen maximal independent set in the graph $G - I_0 - \cdots - I_{n-1}$. We claim that $V(G) = \bigsqcup_{n=0}^\Delta I_n$. Indeed, every vertex not in $\bigsqcup_{n=0}^\Delta I_n$ must have a neighbor in each of $I_0$, \ldots, $I_\Delta$, which is impossible as the maximum degree of $G$ is $\Delta$.
	\end{scproof}

	\section{Proof of Theorem~\ref{theo:cont_LLL}}\label{sec:proof_LLL}
	
	\subsection{First observations}


	Call a CSP $\B$ \emphd{bounded} if $\vdeg(\B)$ and $\ord(\B)$ are both finite. Given a CSP $\B \colon X \to^? k$, define a graph $G_\B$ with vertex set $X$ by making two distinct vertices $x$, $y$ adjacent if and only if there is a constraint $B \in \B$ such that $\set{x,y} \subseteq \dom(B)$.  
	
	\begin{lemma}\label{lemma:cont_graph}
		If $\B \colon X \to^? k$ is a bounded continuous CSP on a zero-dimensional Polish space $X$, then the graph $G_\B$ is continuous.
	\end{lemma}
	\begin{scproof}
		Set $G \defeq G_\B$ and let $U \subseteq X$ be a clopen set. A vertex $x_1$ is in $\Nbhd_G(U)$ if and only if there are some $2 \leq i \leq n \leq \ord(\B)$ and $B \subseteq k^{\set{1,\ldots, n}}$ such that:
		\[
			\text{there are $x_2 \in X$, \ldots, $x_i \in U$, \ldots $x_n \in X$ such that $x_1$, \ldots, $x_n$ are distinct and $B(x_1, \ldots, x_n) \in \B$}.
		\]
		This shows that $\Nbhd_G(U)$ is a union of finitely many clopen sets, hence it is itself clopen.
	\end{scproof}

	
	Let $X$ be a set and let $g \colon X' \to k$ be a coloring with domain $X' \subseteq X$. Given an $(X,k)$-constraint $B$ with domain $D$, let $B/g$ be the constraint with domain $\dom(B/g) \defeq D \setminus X'$ given by
	\[
	B/g \,\defeq\, \set{\phi \colon D \setminus X' \to k \,:\, \rest{g}{D \cap X'} \sqcup \phi \,\in\, B}.
	\]
	In other words, $\phi \in B/g$ if and only if the coloring $g \sqcup \phi$ violates $B$. Here it is possible that $D \subseteq X'$, in which case $\dom(B/g) = \0$; more specifically, $B/g = \set{\0}$ if $g$ violates $B$, and $B/g = \0$ otherwise. \ep{Note that the constraint $\set{\0}$ has probability $1/k^0 = 1$ and is violated by every coloring, while the constraint $\0$ has probability $0$ and is always satisfied.} For a CSP $\B \colon X \to^? k$, we define \[\B/g \,\defeq\, \set{B/g \,:\, B \in \B}\] and view $\B/g$ as a CSP on $X \setminus X'$. By construction, $h \colon X \setminus X' \to k$ is a solution to $\B/g$ if and only if $g \sqcup h$ is a solution to $\B$. 
	
	\begin{lemma}\label{lemma:part_cont}
		Let $\B \colon X \to^? k$ be a bounded continuous CSP on a zero-dimensional Polish space $X$. If $X' \subseteq X$ is a clopen set and $g \colon X' \to k$ is continuous, then the CSP $\B/g \colon X \setminus X' \to^? k$ is also continuous.
	\end{lemma}
	\begin{scproof}
		The proof is very similar to the proof of Lemma~\ref{lemma:cont_graph}. Given a set $B \subseteq k^{\set{1, \ldots,n}}$ and clopen subsets $U_2$, \ldots, $U_n \subseteq X \setminus X'$, we have to argue that the following set is clopen:
		\[
		\set{x_1 \in X \setminus X' \,:\, \text{$\exists$ $x_2 \in U_2$, \ldots, $x_n \in U_n$ such that $x_1$, \ldots, $x_n$ are distinct and $B(x_1, \ldots, x_n) \in \B/g$}}.
		\]
		To this end, observe that this set can be written as a union of finitely many clopen sets of the form
		\begin{align*}
		(X \setminus X') \,\cap\, \{x_1 \in X \,:\, &\text{$\exists$ $x_2 \in U_2$, \ldots, $x_n \in U_n$, $x_{n+1} \in g^{-1}(\alpha_{n+1})$, \ldots, $x_m \in g^{-1}(\alpha_m)$}\\
		&\text{ such that $x_1$, \ldots, $x_m$ are distinct and $\tilde{B}(x_1,\ldots, x_m) \in \B$}\},
		\end{align*}
		for some $n \leq m \leq \ord(\B)$, colors $0 \leq \alpha_{n+1}$, \ldots, $\alpha_m \leq k- 1$, and an appropriate set $\tilde{B} \subseteq k^{\set{1,\ldots, m}}$.
	\end{scproof}

	\subsection{Good CSPs and conditional probabilities}

	Call a CSP $\B$ \emphd{good} if it is bounded and for all $B \in \B$,
	\begin{equation}\label{eq:good}
	\P[B] \cdot \vdeg(\B)^{|\dom(B)|} \,<\, 1.
	\end{equation}
	If $\vdeg(\B) = |\dom(B)| = 0$, we interpret the expression $0^0$ appearing in \eqref{eq:good} as $1$. Note that every CSP satisfying \eqref{eq:bound} is good. The following lemma is the main step in the proof of Theorem~\ref{theo:cont_LLL}:

	\begin{lemma}\label{lemma:step}
		Let $\B \colon X \to^? k$ be a good continuous CSP on a zero-dimensional Polish space $X$, and let $I \subseteq X$ be a clopen independent set in $G_\B$. Then there is a continuous coloring $g \colon I \to k$ such that $\B/g$ is good.
	\end{lemma}
	\begin{scproof}
		For brevity, let $G \defeq G_\B$ and $\vdeg \defeq \vdeg(\B)$. Note that $\vdeg(\B/g) \leq \vdeg$ for every $g \colon I \to k$, so it is enough to argue that there is a continuous coloring $g \colon I \to k$ such that
		\begin{equation}\label{eq:new_bound}
		\P[B/g] \cdot \vdeg^{|\dom(B/g)|} \,<\, 1 \quad \text{for all }B \in \B.
		\end{equation}
		For each $x \in I$, let $\B_x \subseteq \B$ denote the set of all constraints $B$ with $x \in \dom(B)$. Note that $|\B_x| \leq \vdeg$. Since $I$ is independent in $G$, $x$ is the unique element of $I \cap \dom(B)$ for each $B \in \B_x$; in particular, the value $\P[B/g]$ only depends on the color $g(x)$. Specifically, for each $B \in \B_x$ and a color $\alpha$, we define
		\[
		\P[B\,\vert\, x \mapsto \alpha] \,\defeq\, \frac{|\set{\phi \in B \,:\, \phi(x) = \alpha}|}{k^{|\dom(B)| - 1}}.
		\]
		Then for any coloring $g \colon I \to k$, $\P[B/g] = \P[B\,\vert\, x \mapsto g(x)]$. We say that a color $\alpha$ is \emphd{good} for $x$ if
		\[
			\P[B\,\vert\, x \mapsto \alpha] \,\leq\, \P[B] \cdot \vdeg \quad \text{for all } B \in \B_x.
		\]
		
		\begin{claim*}\label{claim:claim}
			For each $x \in I$, there is a good color.
		\end{claim*}
		\begin{stepproof}
			Take any $x \in I$ and notice that for each $B \in \B_x$, $\P[B] = (1/k)\sum_{\alpha = 0}^{k-1} \P[B\,\vert\, x \mapsto \alpha]$. This implies that there are fewer than $k/\vdeg$ colors $\alpha$ such that $\P[B\,\vert\, x \mapsto \alpha] > \P[B] \cdot \vdeg$. Since $|\B_x| \leq \vdeg$, there are fewer than $k$ colors that are not good for $x$, as desired.
		\end{stepproof}
	
		Now we define $g \colon I \to k$ by making $g(x)$ be the minimum color that is good for $x$. Since $\B$ is continuous, it is straightforward to check 
		that $g$ is continuous. It remains to verify that \eqref{eq:new_bound} holds. To this end, take any $B \in \B$. If $I \cap \dom(B) = \0$, then $B/g = B$ and \eqref{eq:new_bound} is satisfied automatically \ep{since $\B$ is good}. Otherwise, $B \in \B_x$ for some \ep{unique} $x \in I$, and we can write
		\begin{align*}
			\P[B/g] \cdot \vdeg^{|\dom(B/g)|} \,&=\, \P[B\,\vert\, x \mapsto g(x)] \cdot \vdeg^{|\dom(B)| - 1} \\
			\big[\text{since $g(x)$ is good for $x$}\big]\quad\quad&\leq\, \P[B] \cdot \vdeg \cdot \vdeg^{|\dom(B)| - 1} \\
			&=\, \P[B] \cdot \vdeg^{|\dom(B)|} \\
			\big[\text{since $\B$ is good}\big]\quad\quad&<\, 1. \qedhere
		\end{align*}
	\end{scproof}

	We are now ready to prove the following strengthening of Theorem~\ref{theo:cont_LLL}:
	
	\begin{theo}
		If $\B \colon X \to^? k$ is a good continuous CSP on a zero-dimensional Polish space $X$, then $\B$ has a continuous solution $f \colon X \to k$.
	\end{theo}
	\begin{scproof}
		The graph $G \defeq G_\B$ has $\Delta(G) \leq \vdeg(\B)(\ord(\B) - 1) < \infty$, so, by Lemmas~\ref{lemma:cont_graph} and \ref{lemma:coloring}, there is a partition $X = I_1 \sqcup \ldots \sqcup I_n$ of $X$ into finitely many clopen sets that are independent in $G$. Thanks to Lemma~\ref{lemma:part_cont}, we may iteratively apply Lemma~\ref{lemma:step} to produce a sequence of continuous colorings $g_i \colon I_i \to k$ such that for all $i \leq n$, the CSP $\B/(g_1 \sqcup \ldots \sqcup g_i)$ is good. We claim that $f \defeq g_1 \sqcup \ldots \sqcup g_n$ is a solution to $\B$, as desired. Indeed, suppose $f$ violates a constraint $B \in \B$. Then we have $B/f = \set{\0}$, but this means that $\P[B/f] = 1$, contradicting the fact that the CSP $\B/f$ is good.
	\end{scproof}

	\section{Proofs of Theorems~\ref{theo:STD} and \ref{theo:top_AW}}\label{sec:main_proof}
	
	\subsection{The main lemma}
	
	Recall that $\G$ is a countably infinite group with identity element $\mathbf{1}$. Given an action $\G \acts X$ and a set $S \subset \G$, a subset $A \subseteq X$ is \emphd{$S$-syndetic} if $S^{-1} \cdot A = X$ and \emphd{$S$-separated} if for all distinct $x$, $y \in A$, $y \not \in S \cdot x$. Note that a set $A \subseteq X$ is $S$-separated if and only if it is independent in the Schreier graph $G(X, S)$. If $X$ is a free zero-dimensional Polish $\G$-space, then the neighborhood of a clopen set $U \subseteq X$ in $G(X,S)$ is $((S \cup S^{-1}) \setminus \set{\mathbf{1}}) \cdot U$, which is also clopen. Hence, in this situation the graph $G(X,S)$ is continuous, so we may apply the results of \S\ref{sec:prelim} to it.
	
	Let $\G\acts X$ be an action and let $f \colon X \rightharpoonup k$ be a partial coloring. Given a subset $S \subseteq \G$, we say that two points $x$, $y \in X$ are \emphd{$S$-similar} in $f$, in symbols $x \equiv^S_f y$, if
	\[
	\forall \sigma \in S, \quad \set{\sigma \cdot x,\, \sigma \cdot y} \, \subseteq \, \dom(f) \quad \Longrightarrow \quad f(\sigma \cdot x) \,=\, f(\sigma \cdot y).
	\]
	
	\begin{lemma}\label{lemma:syndetic}
		For every finite set $F \subset \G$, there is a finite set $S \subset \G$ with the following property:
		\noindent Let $X$ be a free zero-dimensional Polish $\G$-space and let $X \defeq C_0 \sqcup C \sqcup U$ be a partition of $X$ into clopen sets such that $C$ is $F$-syndetic and $U$ is $S$-separated. Then, given an element $\mathbf{1} \neq \gamma \in \G$,
		every continuous $2$-coloring $f_0 \colon C_0 \to 2$ can be extended to a continuous $2$-coloring $f \colon C_0 \sqcup C \to 2$ such that 
		\begin{equation}\label{eq:want1}
		\forall x \in X, \quad x \,\not\equiv^S_f \, \gamma \cdot x. 
		\end{equation}
	\end{lemma}

	In the notation of Lemma~\ref{lemma:syndetic}, the set $C_0$ is already \emph{colored}, the set $C$ is the one we need to \emph{color}, and the set $U$ will be left \emph{uncolored}. Lemma~\ref{lemma:syndetic} is analogous to \cite[Lemma~3.9]{STD} and is used in much the same inductive fashion in our proof of Theorem~\ref{theo:STD}. The main novelty of our approach is in the proof of Lemma~\ref{lemma:syndetic}, which uses Theorem~\ref{theo:cont_LLL}.
	
	\begin{scproof}
	Let $F \subset \G$ be a finite set. 
	We may assume that $F$ is symmetric \ep{i.e., $F^{-1} = F$} and $\mathbf{1} \in F$. Let $M$ be any finite symmetric subset of $\G$ with $\mathbf{1} \in M$ of size $|M| = m|F|$, where $m > 0$ is so large that
	\begin{equation}\label{eq:m}
		2^{m} \,>\, (2m|F|)^{500}.
	\end{equation}
	This inequality will only be used on the very last step of the argument, where it will be invoked to ensure that the numerical requirements of Theorem~\ref{theo:cont_LLL} are fulfilled. Let $N \defeq FM \cup MF$. We claim that the conclusion of Lemma~\ref{lemma:syndetic}
	holds for $S \defeq N^5F$.
	
	Let $X$ be a free zero-dimensional Polish $\G$-space and let $X \defeq C_0 \sqcup C \sqcup U$ be a partition of $X$ into clopen sets such that $C$ is $F$-syndetic and $U$ is $S$-separated. Fix a group element $\gamma \neq \mathbf{1}$ and let $\Delta \defeq N^4 F \gamma F N^4 \setminus \set{\mathbf{1}}$.
	By Lemma~\ref{lemma:max}, there is a clopen maximal $N^4$-separated subset $Z$ of $C$. 
	Since $N$ is symmetric and contains $\mathbf{1}$, the maximality of $Z$ means that $C \subseteq N^4 \cdot Z$. Since $C$ is $F$-syndetic, this implies that $Z$ is $N^4F$-syndetic. 
	Let $g \colon C_0 \sqcup (C \setminus (N \cdot Z)) \to 2$ be an arbitrary continuous extension of $f_0$ \ep{for instance, we can set $g(x) \defeq 0$ for all $x \in C \setminus (N \cdot Z)$}. We shall extend $g$ to a continuous coloring $f \colon C_0 \sqcup C \to 2$ such that 
	\begin{equation}\label{eq:want}
	\forall z \in Z \, \forall \delta \in \Delta, \quad z \,\not\equiv^N_f \, \delta \cdot z. 
	\end{equation}
	
	\begin{claim*}
		If $f$ satisfies \eqref{eq:want}, then it also satisfies \eqref{eq:want1}. 
	\end{claim*}
	\begin{stepproof}
		Take any $x \in X$. Since $Z$ is $N^4F$-syndetic, there is $\beta \in N^4 F$ such that $\beta \cdot x \in Z$. Applying \eqref{eq:want} with $z = \beta \cdot x$ and $\delta = \beta \gamma \beta^{-1}$, we get $\beta \cdot x \not\equiv^N_f \beta \gamma \cdot x$. Since $N\beta \subseteq S$, this yields 
		$x \not \equiv^S_f \gamma \cdot x$, as desired. 
	\end{stepproof}

	Extensions of $g$ to $C_0 \sqcup C$ can be encoded by $2^{|N|}$-colorings of $Z$, as follows. A natural number less than $2^{|N|}$ can be identified with a binary sequence of length $|N|$, so a $2^{|N|}$-coloring $h \colon Z \to 2^{|N|}$ can be viewed as an $|N|$-tuple of $2$-colorings $h_1$, \ldots, $h_{|N|} \colon Z \to 2$. Let $N = \set{\nu_1, \ldots, \nu_{|N|}}$ be an enumeration of $N$. Since $X$ is free and $Z$ is $N^4$-separated, each point $x \in N \cdot Z$ can be expressed uniquely as $x = \nu_i \cdot z$ for some $z \in Z$ and $1 \leq i \leq |N|$. Thus, given $h \colon Z \to 2^{|N|}$, we can define $f^h \colon C_0 \sqcup C \to 2$ by the formula
	\[
		f^h(x) \,\defeq\, \begin{cases}
			g(x) &\text{if } x \in C_0 \sqcup (C \setminus (N \cdot Z));\\
			h_i(z) &\text{if } x \in C \text{ and } x = \nu_i \cdot z \text{ for } z \in Z \text{ and } 1 \leq i \leq |N|.
		\end{cases}
	\]
	In other words, for each $z \in Z$, the color $h(z) \in 2^{|N|}$ encodes the restriction of $f^h$ to the set $C \cap (N \cdot z)$. This encoding is generally not one-to-one: unless $N \cdot z \subseteq C$, the sequence $h_1(z)$, \ldots, $h_{|N|}(z)$ includes some redundant bits. Nevertheless, choosing $h(z)$ uniformly at random does correspond to picking 
	the restriction of $f^h$ to $C \cap (N \cdot z)$
	uniformly at random form the set of all $2$-colorings $C \cap (N \cdot z) \to 2$. Notice also that if $h$ is continuous, then so is $f^h$.
	
	To apply Theorem~\ref{theo:cont_LLL}, we now need to define a constraint satisfaction problem $\B \colon Z \to^? 2^{|N|}$ such that $h \colon Z \to 2^{|N|}$ is a solution to $\B$ if and only if $f^h$ satisfies \eqref{eq:want}, i.e.,
	\[
		\text{$h$ is a solution to $\B$} \quad \Longleftrightarrow \quad \forall z \in Z \, \forall \delta \in \Delta, \quad z \,\not\equiv^N_{f^h} \, \delta \cdot z.
	\]
	To this end, observe that the truth of the statement $z \not\equiv^N_{f^h} \delta \cdot z$ only depends on the restriction of $f^h$ to $(N \cdot z) \cup (N\delta \cdot z)$. Thus, for each $z \in Z$ and $\delta \in \Delta$, there is a constraint $B_{z, \delta}$ with domain
	\[
		\dom(B_{z,\delta}) \,\defeq\, \set{z' \in Z \,:\, (N \cdot z') \cap ((N \cdot z) \cup (N\delta \cdot z)) \neq \0} \,=\, Z \cap \left((N^2 \cup N^2\delta) \cdot z\right)
	\]
	such that $h$ satisfies $B_{z, \delta}$ if and only if $z \not\equiv^N_{f^h} \delta\cdot z$. We then let $\B\defeq \set{B_{z,\delta} \,:\, z \in Z \text{ and } \delta \in \Delta}$. It is clear from the definition that the CSP $\B$ is continuous.
	
	\begin{claim*}
		$\ord(\B) \leq 2$.
	\end{claim*}
	\begin{stepproof}
		Since $Z$ is $N^4$-separated, $|Z \cap (N^2 \cdot x)| \leq 1$ for all $x \in X$. Hence, for any $z \in Z$ and $\delta \in \Delta$, there are at most $2$ elements in $Z \cap \left((N^2 \cup N^2\delta) \cdot z\right)$, i.e., $|\dom(B_{z,\delta})| \leq 2$, as desired.
	\end{stepproof}

	\begin{claim*}
		$\vdeg(\B) \leq 2^{11}m^{10}|F|^{22}$.
	\end{claim*}
	\begin{stepproof}
		Take any $z' \in Z$. We need to bound the number of pairs $(z, \delta) \in Z \times \Delta$ such that $z' \in \dom(B_{z,\delta})$. Recall that $N = FM \cup MF$, where $|M| = m|F|$, so $|N| \leq 2m|F|^2$. Hence, $|\Delta| \leq |N^4F\gamma F N^4| \leq 2^8 m^8|F|^{18}$. Once $\delta$ is fixed, $z$ must satisfy $z' \in (N^2 \cup N^2 \delta) \cdot z$, i.e., $z \in (N^2 \cup \delta^{-1}N^2) \cdot z'$, so there are at most $8m^2|F|^4$ such $z$. Thus, the number of choices for $(z, \delta)$ is at most $2^8m^8|F|^{18} \cdot 8m^2|F|^4 = 2^{11}m^{10}|F|^{22}$. 
	\end{stepproof}
	
	\begin{claim*}
		$\pr(\B) \leq 2^{-m/6}$.
	\end{claim*}
	\begin{stepproof}
		Take any $z \in Z$ and $\delta \in \Delta$. For brevity, let $y \defeq \delta \cdot z$. We need to show that $\P[B_{z,\delta}] \leq 2^{-m/6}$, i.e., the probability that $z$ is $N$-similar to $y$ in a random extension $f$ of $g$ to $C_0 \sqcup C$ is at most $2^{-m/6}$.
		
		Call an element $\nu \in N$ \emphd{eligible} if $\nu \cdot z \in C$ and $\nu \cdot y \in C_0 \sqcup C$. Let $E$ be the set of all eligible $\nu \in N$. Note that if $\nu$ is eligible, then $\nu \cdot z$ is uncolored in $g$ but becomes colored in $f$, and $\nu \cdot y$ is also colored in $f$ \ep{but it may or may not be already colored in $g$}. The color $f(\nu \cdot z)$ is chosen randomly, so the probability that $f(\nu \cdot z) = f(\nu \cdot y)$ is exactly $1/2$, regardless of whether $\nu \cdot y$ is already colored in $g$. 
		
		Since $C$ is $F$-syndetic and $N \supseteq FM$, we have $|C \cap (N \cdot z)| \geq |M|/|F| = m$, and since $U$ is $S$-separated and $S \supseteq N^2$, $|(N \cdot y) \cap U| \leq 1$. Therefore, $|E| \geq m - 1 \geq m/2$. Let $G$ be the graph with vertex set $(N \cdot z) \cup (N \cdot y)$ in which we put an edge between $\nu \cdot z$ and $\nu \cdot y$ for each $\nu \in E$. The maximum degree of $G$ is at most $2$, so we can pick a subset $E' \subseteq E$ of size $|E'| \geq |E|/3 \geq m/6$ such that the pairs $\set{\nu \cdot z, \nu \cdot y}$, $\nu \in E'$, are pairwise disjoint. When $f$ is chosen randomly, the events $f(\nu \cdot z) = f(\nu \cdot y)$ for distinct $\nu \in E'$ are mutually independent, so the probability that they all occur simultaneously is $2^{-|E'|} \leq 2^{-m/6}$, which gives us the desired upper bound on the probability that $z$ is $N$-similar to $y$ in $f$.
	\end{stepproof}
	
	And now we are done: by Theorem \ref{theo:cont_LLL}, $\B$ has a continuous solution as long as
	\[
		\pr(\B) \cdot \vdeg(\B)^{\ord(\B)} \,\leq\, 2^{-m/6} \cdot (2^{11}m^{10}|F|^{22})^2 \,=\, 2^{-m/6} \cdot 2^{22}m^{20}|F|^{44}\,<\, 1,
	\]
	which holds by \eqref{eq:m}.
	\end{scproof}


	\subsection{Proof of Theorem~\ref{theo:STD}}\label{subsec:proof_ST-D}
	
	For the reader's convenience, we state Theorem~\ref{theo:STD} again:
	
	\begin{theocopy}{theo:STD}
		If $\G \acts X$ is a free Borel action of $\G$ on a standard Borel space $X$, then there is a $\G$-equivariant Borel map $\pi \colon X \to Y$, where $Y \subset 2^\G$ is a free subshift.
	\end{theocopy}
	
	To prove Theorem~\ref{theo:STD}, we shall first define a free subshift $Y \subset 2^\G$ and then iteratively apply Lemma~\ref{lemma:syndetic} to construct a desired $\G$-equivariant Borel map $\pi \colon X \to Y$.
	
	We start by recursively defining a sequence of finite sets $H_0$, $F_0$, $S_0$, $H_1$, $F_1$, $S_1$, \ldots{} $\subset \G$ as follows. Let $H_0$ be an arbitrary nonempty finite subset of $\G$. Once $H_n$ is defined, let $\delta_n$ be any group element such that $H_n \cap (H_n\delta_n) = \0$ \ep{such $\delta_n$ exists since $\G$ is infinite} and set $F_n \defeq H_n \cup (H_n\delta_n)$. Next, let $S_n$ be the set $S$ produced by Lemma~\ref{lemma:syndetic} applied with $F = F_n$. Upon replacing $S_n$ with a superset if necessary, we may additionally assume that $S_n$ is symmetric and $S_n \supseteq F_n^{-1}F_n$. Finally, we let $H_{n+1} \defeq S_nH_n$. The following claim explains why the sets $H_n$, $F_n$, and $S_n$ are defined in this manner.
	
	\begin{big_claim}\label{claim:split}
		Let $X$ be a free zero-dimensional Polish $\G$-space and let $W \subseteq X$ by an $H_n$-syndetic clopen set. Then there is a partition $W = C \sqcup U$ into two clopen sets such that:
		\begin{itemize}
			\item the set $C$ is $F_n$-syndetic;
			\item the set $U$ is $S_n$-separated and $H_{n+1}$-syndetic.
		\end{itemize}
	\end{big_claim}
	\begin{scproof}
		By Lemma~\ref{lemma:max}, we can let $U$ be a clopen maximal $S_n$-separated subset of $W$ and define $C \defeq W \setminus U$. Since $S_n$ is symmetric and contains $\mathbf{1}$, the maximality of $U$ means that $W \subseteq S_n \cdot U$, and since $W$ is $H_n$-syndetic and $H_{n+1} = S_n H_n$, this implies that $U$ is $H_{n+1}$-syndetic, as claimed.
		
		To see that $C$ is $F_n$-syndetic, take any $x \in X$. We need to argue that $F_n \cdot x$ contains a point in $C$.  Recall that $F_n = H_n \cup (H_n \delta_n)$. Since $W$ is $H_n$-syndetic, the sets $H_n \cdot x$ and $H_n \delta_n \cdot x$ each contain a point in $W$. Since the sets $H_n$ and $H_n \delta_n$ are disjoint, we have $|(F_n \cdot x) \cap W| \geq 2$. On the other hand, $|(F_n \cdot x) \cap U| \leq 1$ since $U$ is $F_n^{-1}F_n$-separated. Therefore, $|(F_n \cdot x) \cap C| \geq 1$, as desired.
	\end{scproof}

	Fix an arbitrary enumeration $\gamma_0$, $\gamma_1$, \ldots{} of the non-identity elements of $\G$. For each $n \in \N$, let $Y_n \subset 2^\G$ be the set of all $2$-colorings $y \colon \G \to 2$ such that
	\[
		\exists \sigma \in S_n \text{ with } y(\sigma) \,\neq\, y(\sigma\gamma_n).
	\]
	The set $Y_n$ is clopen, and if $y \in Y_n$, then $\gamma_n \cdot y \neq y$. Hence, the set $Y' \defeq \bigcap_{n = 0}^\infty Y_n$ is closed and every point $y \in Y'$ has trivial stabilizer. Finally, we define $Y \defeq \bigcap_{\delta \in \G} (\delta \cdot Y')$. The set $Y$ is closed, $\G$-invariant, and contained in $Y' \subseteq \Free(2^\G)$, so $Y$ is a free subshift \ep{although we have not yet shown that $Y$ is nonempty}.
	
	Now let $\G \acts X$ be a free Borel action of $\G$ on a standard Borel space $X$. It follows from standard results in descriptive set theory that there is a compatible zero-dimensional Polish topology $\tau$ on $X$ with respect to which the action $\G \acts X$ is continuous \cite[\S13]{KechrisDST}. Iterative applications of Claim~\ref{claim:split} yield a sequence of clopen subsets $U_0$, $C_0$, $U_1$, $C_1$, \ldots{} of $X$ such that $U_0 = X$ and for all $n \in \N$,
	\begin{itemize}
		\item $U_n = C_n \sqcup U_{n+1}$; and
		\item the set $C_n$ is $F_n$-syndetic, while $U_{n+1}$ is $S_n$-separated and $H_{n+1}$-syndetic.
	\end{itemize}
	Next we use Lemma~\ref{lemma:syndetic} repeatedly to obtain an increasing sequence $f_0 \subseteq f_1 \subseteq \ldots$ such that for each $n \in \N$, $f_n \colon C_0 \sqcup \ldots \sqcup C_n \to 2$ is a continuous $2$-coloring satisfying
	\begin{equation}\label{eq:not_sim}
		\forall x \in X, \quad x \,\not\equiv^{S_n}_{f_n} \, \gamma_n \cdot x.
	\end{equation}
	Let $f \colon X \to 2$ be an arbitrary Borel extension of $\bigcup_{n=0}^\infty f_n$ \ep{e.g., we may set $f(x) \defeq 0$ for all $x \not \in \bigsqcup_{n=0}^\infty C_n$}. Define a $\G$-equivariant Borel map $\pi_f \colon X \to 2^\G$ by setting $\pi_f(x)(\gamma) \defeq f(\gamma \cdot x)$ for all $x \in X$ and $\gamma \in \G$. We claim that $\pi_f(x) \in Y$ for all $x \in X$, as desired. Indeed, since $\pi_f$ is $\G$-equivariant, it suffices to argue that $\pi_f(x) \in Y_n$ for all $x \in X$ and $n \in \N$, i.e., that for all $x \in X$ and $n \in \N$,
	\[
		\exists \sigma \in S_n \text{ with } \pi_f(x)(\sigma) \,\neq\, \pi_f(x)(\sigma \gamma_n).
	\]
	Using the definition of $\pi_f$, we can rewrite the latter statement as
	\[
	\exists \sigma \in S_n \text{ with } f(\sigma \cdot x) \,\neq\, f(\sigma \gamma_n \cdot x),
	\]
	which holds by \eqref{eq:not_sim} since $f$ is an extension of $f_n$.
	
	\subsection{Proof of Theorem~\ref{theo:top_AW}}\label{subsec:proof_top_AW}

	Let us state Theorem~\ref{theo:top_AW} again:
	
	\begin{theocopy}{theo:top_AW}
		If $X$ is a nonempty free zero-dimensional Polish $\G$-space, then $\Free(2^\G) \preccurlyeq X$.
		
		\smallskip
		
		\noindent Explicitly, given any $k \in \N^+$, a finite subset $F \subset \G$, and a continuous $k$-coloring $f \colon \Free(2^\G) \to k$, there is a continuous $k$-coloring $g \colon X \to k$ such that $\Pat_F(X, g) = \Pat_F(\Free(2^\G), f)$.
	\end{theocopy}

	Our proof of Theorem~\ref{theo:top_AW} is a modification of the proof of Theorem~\ref{theo:STD} presented in \S\ref{subsec:proof_ST-D}. To begin with, fix $k \in \N^+$, a finite subset $F \subset \G$, and a continuous $k$-coloring $f \colon \Free(2^\G) \to k$. The following clopen sets from a base for the topology on $2^\G$: 
	\[
		U(s) \,\defeq\, \set{x \in 2^\G \,:\, x(\gamma) = s(\gamma) \text{ for all } \gamma \in \dom(s)},
	\]
	where $s$ is a $2$-pattern \ep{i.e., a partial mapping $s \colon \G \rightharpoonup 2$ whose domain is a finite subset of $\G$}. Given a finite set $D \subset \G$ and a point $x \in \Free(2^\G)$, we say that $D$ \emphd{$f$-determines} $x$ if for all $z \in \Free(2^\G)$,
	\[
		\forall \delta \in D, \, z(\delta) = x(\delta) \qquad \Longrightarrow \qquad f(z) = f(x).\qedhere
	\]
	The continuity of $f$ is then equivalent to the following assertion:
	\begin{big_claim}\label{claim:cont}
		For each $x \in \Free(2^\G)$, there is a finite set $D \subset \G$ that $f$-determines $x$. \qed
	\end{big_claim}
	
	
	\begin{big_claim}\label{claim:sp}
		For each $k$-pattern $p \in \Pat_F(\Free(2^\G),f)$, there is a $2$-pattern $s_p$ such that for all $z \in \Free(2^\G)$,
		\[
			z \in U(s_p) \qquad \Longrightarrow \qquad \forall \gamma \in F, \, f(\gamma \cdot z) = p(\gamma).
		\]
	\end{big_claim}
	\begin{scproof}
		Since $p$ occurs in $f$, there is some $x \in \Free(2^\G)$ such that $f(\gamma \cdot x) = p(\gamma)$ for all $\gamma \in F$. Claim~\ref{claim:cont} yields a finite set $D$ such that for all $z \in \Free(2^\G)$,
		\[
		\forall \delta \in D, \, z(\delta) = x(\delta) \qquad \Longrightarrow \qquad \forall \gamma \in F, \, f(\gamma \cdot z) = p(\gamma).
		\]
		Thus, we may take $s_p$ be the $2$-pattern with domain $D$ given by $s_p(\delta) \defeq x(\delta)$ for all $\delta \in D$.
	\end{scproof}

	Set $D \defeq \bigcup \set{\dom(s_p) \,:\, p \in \Pat_F(\Free(2^\G), f)}$ \ep{where $s_p$ is the $2$-pattern given by Claim~\ref{claim:sp}} and let $H_0$ be an arbitrary symmetric finite subset of $\G$ with $|H_0| > |D|$. Next we recursively build a sequence of finite sets $H_0$, $F_0$, $S_0$, $H_1$, $F_1$, $S_1$, \ldots{} $\subset \G$ in the same way we did in \S\ref{subsec:proof_ST-D}. That is, once $H_n$ is defined, we let $\delta_n$ be any group element such that $H_n \cap (H_n\delta_n) = \0$ and set $F_n \defeq H_n \cup (H_n\delta_n)$. Then we let $S_n$ be the set $S$ produced by Lemma~\ref{lemma:syndetic} applied with $F = F_n$. Upon replacing $S_n$ with a superset if necessary, we may additionally assume that $S_n$ is symmetric and $S_n \supseteq F_n^{-1}F_n$. Finally, we let $H_{n+1} \defeq S_nH_n$. The following is a restatement of Claim~\ref{claim:split}:
	
	\begin{big_claim}\label{claim:split1}
		Let $X$ be a free zero-dimensional Polish $\G$-space and let $W \subseteq X$ by an $H_n$-syndetic clopen set. Then there is a partition $W = C \sqcup U$ into two clopen sets such that:
		\begin{itemize}
			\item the set $C$ is $F_n$-syndetic;
			\item the set $U$ is $S_n$-separated and $H_{n+1}$-syndetic.
		\end{itemize}
	\end{big_claim}
	\begin{scproof}
		See the proof of Claim~\ref{claim:split} in \S\ref{subsec:proof_ST-D}.
	\end{scproof}

	As in \S\ref{subsec:proof_ST-D}, we now fix an arbitrary enumeration $\gamma_0$, $\gamma_1$, \ldots{} of the non-identity elements of $\G$. For each $n \in \N$, let $Y_n \subset 2^\G$ be the set of all $2$-colorings $y \colon \G \to 2$ such that
	\[
	\exists \sigma \in S_n \text{ with } y(\sigma) \,\neq\, y(\sigma\gamma_n).
	\]
	Let $Y \defeq \bigcap_{n=0}^\infty\bigcap_{\delta \in \G} (\delta \cdot Y_n)$. As discussed in \S\ref{subsec:proof_ST-D}, $Y$ is a free subshift \ep{it is also shown there that $Y$ is nonempty}. For each $N \in \N$, we also define
	\[
		Y_{\leq N} \,\defeq\, \bigcap_{n = 0}^N\bigcap_{\delta \in \G} (\delta \cdot Y_n).
	\]
	Then $Y_{\leq N}$ is a subshift and $Y=\bigcap_{N=0}^\infty Y_{\leq N}$, where the intersection is decreasing. Note that $Y_{\leq N}$ need not be free; in particular, $f$ may not be defined on all of $Y_{\leq N}$. Nevertheless, for large enough $N$, it is possible to define a continuous $k$-coloring $f^\ast \colon Y_{\leq N} \to k$ that, in some sense, approximates $f$:
	
	\begin{big_claim}\label{claim:approx}
		There exist $N \in \N$ and a continuous $k$-coloring $f^\ast \colon Y_{\leq N} \to k$ such that for each $z \in Y_{\leq N}$, there is $y \in Y$ with the following properties:
		\begin{itemize}
			\item for all $\delta \in D$, $z(\delta) = y(\delta)$; and
			\item for all $\gamma \in F$, $f^\ast(\gamma \cdot z) = f(\gamma \cdot y)$.
		\end{itemize}
	\end{big_claim}
	\begin{scproof}
		First we argue that there is a finite set $L \subset \G$ that $f$-determines every point $y \in Y$.	For each finite set $L \subset \G$, let $V_L$ be the set of all points $y \in Y$ that are $f$-determined by $L$. Each set $V_L$ is relatively open in $Y$. Moreover, by Claim~\ref{claim:cont}, the union of all the sets $V_L$ is $Y$. Since $Y$ is compact, this implies that there is a finite collection $L_1$, \ldots, $L_r$ of finite subsets of $\G$ such that $Y = V_{L_1} \cup \ldots \cup V_{L_r}$. Then every point $y \in Y$ is $f$-determined by $L \defeq L_1 \cup \ldots \cup L_r$, as desired.
		
		Next we observe that there is $N \in \N$ such that for each $z \in Y_{\leq N}$,
		\begin{equation}\label{eq:approx}
		\exists y \in Y \text{ such that } \forall \delta \in D \cup L \cup  LF, \, z(\delta) = y(\delta).
		\end{equation}
		Indeed, let $Q$ be the set of all $z \in 2^\G$ for which \eqref{eq:approx} fails. Then $Q$ is a clopen subset of $2^\G$ and $Q \cap Y = \0$. Since $2^\G$ is compact and $Y = \bigcap_{N=0}^\infty Y_{\leq N}$, there must exist some $N \in \N$ with $Q \cap Y_{\leq N} = \0$, as desired.
		
		Finally, we define a $k$-coloring $f^\ast \colon Y_{\leq N} \to k$ as follows:
		\begin{align*}
			f^\ast(z) = c \quad \vcentcolon&\Longleftrightarrow \quad \exists y \in Y \text{ such that } f(y) =c \text{ and } \forall \delta \in L, \, z(\delta) = y(\delta)\\
			&\Longleftrightarrow \quad \forall y \in Y, \text{ we have } \left(\forall \delta \in L, \, z(\delta) = y(\delta)\right) \ \Longrightarrow \ f(y) = c.
		\end{align*}
		The two definitions given above are equivalent since every $y \in Y$ is $f$-determined by $L$. By construction, $L$ also $f^\ast$-determines every $z \in Y_{\leq N}$, so $f^\ast$ is continuous. Now consider any $z \in Y_{\leq N}$. By \eqref{eq:approx}, there is $y \in Y$ such that for all $\delta \in D \cup L \cup LF$, $z(\delta) = y(\delta)$, and it is clear that $y$ has the desired properties.
	\end{scproof}

	Now let $X$ be a nonempty free zero-dimensional Polish $\G$-space. Fix $N \in \N$ and $f^\ast \colon Y_{\leq N} \to k$ given by Claim~\ref{claim:approx}. We shall construct a continuous $k$-coloring $g \colon X \to k$ such that $\Pat_F(X, g) = \Pat_F(\Free(2^\G), f)$ by first building a continuous $\G$-equivariant map $\pi \colon X \to Y_{\leq N}$ and then setting $g \defeq f^\ast \circ \pi$.
	
	We start our construction by letting $W \subseteq X$ be a clopen maximal $D^{-1}H_0^2D$-separated subset of $X$ \ep{which exists by Lemma~\ref{lemma:max}}. Since $X$ is free and nonempty, every $\G$-orbit in $X$ intersects $W$ in infinitely many points, so $W$ is infinite. Thus, we may partition $W$ as $W = \bigsqcup_p W_p$, where the union is over all $p \in \Pat_F(\Free(2^\G), f)$ and each $W_p$ is nonempty and clopen. Let $B_p \defeq \dom(s_p) \cdot W_p$ and $B \defeq \bigsqcup_p B_p$ \ep{the union is disjoint since $W$ is $D^{-1}D$-separated} and define a continuous $2$-coloring $b \colon B \to 2$ by
	\begin{equation}\label{eq:b}
		b(\delta \cdot w) \,\defeq\, s_p(\delta) \text{ for all $p \in \Pat_F(\Free(2^\G), f)$, $w \in W_p$, and $\delta \in \dom(s_p)$}. 
	\end{equation}
	Property \eqref{eq:b} will be eventually used to show that $\Pat_F(X, g) \supseteq \Pat_F(\Free(2^\G), f)$.
	
	To continue our construction, 
	we need to make sure that $X \setminus B$ is syndetic:
	
	\begin{big_claim}\label{claim:base}
		The set $X \setminus B$ is $H_0$-syndetic.
	\end{big_claim}
	\begin{scproof}
		Take any $x \in X$. Since $W$ is $D^{-1}H_0^2D$-separated, there is at most one $w \in W$ such that $(D \cdot w) \cap (H_0 \cdot x) \neq \0$, so $|(D \cdot W) \cap (H_0 \cdot x)| \leq |D| < |H_0|$. Since $B \subseteq D \cdot W$, this implies $(H_0 \cdot x) \setminus B \neq \0$.
	\end{scproof}

	Claim~\ref{claim:base} allows us to iteratively apply Claim~\ref{claim:split1} in order to obtain a sequence of clopen subsets $U_0$, $C_0$, $U_1$, $C_1$, \ldots{} of $X$ such that $U_0 = X \setminus B$ and for all $n \in \N$,
	\begin{itemize}
		\item $U_n = C_n \sqcup U_{n+1}$; and
		\item the set $C_n$ is $F_n$-syndetic, while $U_{n+1}$ is $S_n$-separated and $H_{n+1}$-syndetic.
	\end{itemize}
	We can then use Lemma~\ref{lemma:syndetic} repeatedly to obtain an increasing sequence $b \subseteq h_0 \subseteq h_1 \subseteq \ldots$ such that for each $n \in \N$, $h_n \colon B \sqcup C_0 \sqcup \ldots \sqcup C_n \to 2$ is a continuous $2$-coloring satisfying
	\begin{equation}\label{eq:not_sim1}
	\forall x \in X, \quad x \,\not\equiv^{S_n}_{h_n} \, \gamma_n \cdot x.
	\end{equation}
	
	Recall that $N \in \N$ and $f^\ast \colon Y_{\leq N} \to k$ are given by Claim~\ref{claim:approx}. Let $h \colon X \to 2$ be an arbitrary continuous extension of $h_N$ \ep{e.g., we may set $h(x) \defeq 0$ for all $x \not \in \dom(h_N)$} and define a $\G$-equivariant continuous map $\pi_h \colon X \to 2^\G$ by setting $\pi_h(x)(\gamma) \defeq h(\gamma \cdot x)$ for all $x \in X$ and $\gamma \in \G$. Condition~\eqref{eq:not_sim1} ensures that $\pi_h(x) \in Y_{\leq N}$ for all $x \in X$, so we can define a continuous $k$-coloring $g \colon X \to k$ via $g \defeq f^\ast \circ \pi_h$.
	
	\begin{big_claim}\label{claim:sup}
		$\Pat_F(X, g) \supseteq \Pat_F(\Free(2^\G), f)$.
	\end{big_claim}
	\begin{scproof}
		Consider any $p \in \Pat_F(\Free(2^\G), f)$. Take an arbitrary point $w \in W_p$ and let $z \defeq \pi_h(w) \in Y_{\leq N}$. Note that for all $\gamma \in \G$, $g(\gamma \cdot w) = f^\ast(\gamma \cdot z)$. By Claim~\ref{claim:approx}, there is $y \in Y$ such that:
		\begin{enumerate}[label=\ep{\itshape\alph*}]
			\item\label{item:a} for all $\delta \in D$, $z(\delta) = y(\delta)$; and
			\item\label{item:b} for all $\gamma \in F$, $f^\ast(\gamma \cdot z) = f(\gamma \cdot y)$.
		\end{enumerate}
		By \eqref{eq:b}, since $h$ extends $b$, we have $z(\delta) = h(\delta \cdot w) = b(\delta \cdot w) = s_p(\delta)$ for all $\delta \in \dom(s_p)$, i.e., $z \in U(s_p)$. By \ref{item:a}, $y \in U(s_p)$ as well, so for all $\gamma \in F$,
		 \[
		 	g(\gamma \cdot w) \,=\, f^\ast(\gamma \cdot z) \,=\, f(\gamma \cdot y) \,=\, p(\gamma),
		 \]
		 where the second equality holds by \ref{item:b}, and the third by Claim~\ref{claim:sp} and since $y \in U(s_p)$. This shows that $p$ appears in $g$, as desired.
	\end{scproof}

	\begin{big_claim}\label{claim:sub}
		$\Pat_F(X, g) \subseteq \Pat_F(\Free(2^\G), f)$.
	\end{big_claim}
	\begin{scproof}
		Take any $p \in \Pat_F(X, g)$ and let $x \in X$ be such that $g(\gamma \cdot x) = p(\gamma)$ for all $\gamma \in F$. Let $z \defeq \pi_h(x) \in Y_{\leq N}$, so $g(\gamma \cdot x) = f^\ast(\gamma \cdot z)$ for all $\gamma \in \G$. By Claim~\ref{claim:approx}, there is $y \in Y$ such that:
		\begin{itemize}
			\item for all $\gamma \in F$, $f^\ast(\gamma \cdot z) = f(\gamma \cdot y)$.
		\end{itemize}
		Then for all $\gamma \in F$, $f(\gamma \cdot y) = f^\ast(\gamma \cdot z) = g(\gamma \cdot x) = p(\gamma)$, which shows that $p$ appears in $f$, as desired.
	\end{scproof}
	
	Claims~\ref{claim:sup} and \ref{claim:sub} yield $\Pat_F(X, g) = \Pat_F(\Free(2^\G), f)$, and the proof of Theorem~\ref{theo:top_AW} is complete.
	
	
	\section{Combinatorial results}\label{sec:corls}
	
	\subsection{Local colorings of special subshifts}
	
	In this subsection we prove a certain technical result \ep{namely Lemma~\ref{lemma:F_loc_1}} that will be later used to derive Theorems \ref{theo:compact} and \ref{theo:LOCAL}.
	
	Given a subshift $X \subseteq n^\G$, a finite subset $F \subset \G$, and an integer $k \geq 1$, we say that a $k$-coloring $f \colon X \to k$ is \emphd{$F$-local} if for all $x \in X$, the value $f(x)$ is determined by the restriction of $x$ to $F$, i.e., if there is a mapping $\rho \colon n^F \to k$ such that for all $x \in X$, $f(x) = \rho \left((x(\sigma))_{\sigma \in F}\right)$. \ep{In the terminology of \S\ref{subsec:proof_top_AW}, this means that the set $F$ $f$-determines every $x \in X$.} Note that every $F$-local coloring is continuous. Conversely, if $f \colon X \to k$ is continuous, then, due to the compactness of $X$, there is a finite set $F \subset \G$ such that $f$ is $F$-local. 
	
	Let $D$ be a finite subset of $\G$ and let $n \geq 1$ be an integer. Define a subshift $X_{D, n} \subseteq n^\G$ as follows:
	\[
	X_{D, n} \,\defeq\, \set{x \in n^\G \,:\, \text{for all $\gamma \in \G$ and $\sigma \in D \setminus \set{\mathbf{1}}$, we have $x(\gamma) \neq x(\sigma \gamma)$}}.
	\]
	In other words, the elements of $X_{D,n}$ are the proper $n$-colorings of the Cayley graph $G(\G, D)$. 
	\ep{Note that $X_{D,n}$ may be empty if $n$ is too small.} The main result of this subsection allows us to build $F$-local colorings of $X_{D,n}$ with some control over the set $F$:
	
	\begin{lemma}[\textls{Local colorings of $X_{D,n}$}]\label{lemma:F_loc_1}
		Let $\Pat$ be a finite set of $k$-patterns such that every free zero-dimensional Polish $\G$-space admits a continuous $\Pat$-avoiding $k$-coloring. Then there is a finite set $F \subset \G$ with the following property:
		
		\smallskip
		
		Let $n \geq 2$ and let $D \subset \G$ be a finite set such that $F \subseteq D$. Set $F^\ast \defeq F^{\log^\ast n}$. Then the subshift $X_{D,n}$ admits an $F^\ast$-local $\Pat$-avoiding $k$-coloring. 
 	\end{lemma}
 
 	
 	In our proof of Lemma~\ref{lemma:F_loc_1} we shall rely on the following fact, which follows from a construction due to Cole and Vishkin \cite{ColeVishkin}:
 	
 	\begin{lemma}[{\cite[\S3.4]{BE}}]\label{lemma:F_loc_col_1}
 		Let $\gamma \in \G \setminus \set{\mathbf{1}}$ and let $D \subset \G$ be a finite set with $\gamma \in D$. Take $n \geq 2$ and define $F^\ast \defeq \set{\mathbf{1}, \gamma}^{\log^\ast n + 2}$. Then the Schreier graph $G(X_{D,n}, \set{\gamma})$ admits an $F^\ast$-local proper $6$-coloring.
 	\end{lemma}
 	
 	The construction in \cite[\S3.4]{BE} is in the language of distributed algorithms, so, for completeness, we provide its translation into our setting in the \hyperref[app]{appendix}. We remark that it is possible to reduce the number of colors in Lemma~\ref{lemma:F_loc_col_1} from $6$ to $3$ \ep{at the cost of replacing $\log^\ast n + 2$ by $\log^\ast n + C$ for some other constant $C$}, but this will not be needed for our purposes.
 	Since the graph $G(X_{D,n}, \set{\gamma})$ has maximum degree at most $2$, Lemma~\ref{lemma:coloring} already yields a \emph{continuous} proper $3$-coloring of $G(X_{D,n}, \set{\gamma})$. Lemma~\ref{lemma:F_loc_col_1} additionally specifies the set $F^\ast$ such that the resulting coloring is {$F^\ast$-local}.
 
 	\begin{scproof}[Proof of Lemma~\ref{lemma:F_loc_1}]
 		This argument is inspired by Elek's proof of \cite[Theorem 2]{Elek}. Enumerate the non-identity elements of $\G$ as $\gamma_1$, $\gamma_2$, \ldots{} and let $X_i \defeq X_{\set{\gamma_i}, 6}$. Consider the product space $X \defeq \prod_{i=1}^\infty X_i$, equipped with the diagonal action of $\G$. Then $X$ is a compact zero-dimensional Polish $\G$-space. Furthermore, $X$ is free since $\gamma_i \cdot x \neq x$ for all $x \in X_i$. Hence, by the assumptions on $\Pat$, there is a continuous $\Pat$-avoiding $k$-coloring $f \colon X \to k$. The following sets form a base for the topology on $X$:
 		\begin{equation}\label{eq:X_base}
 			\set{x = (x_1, x_2, \ldots) \in X \,:\, x_i(\delta) = s_i(\delta) \text{ for all } 1 \leq i \leq N \text{ and } \delta \in R},
 		\end{equation}
 		where $N$ is a natural number, $R \subset \G$ is a finite set, and $s_1 \colon R \to 6$, \ldots, $s_N \colon R \to 6$ are $6$-patterns. Therefore, each $x \in X$ has a clopen neighborhood of the form \eqref{eq:X_base} on which $f$ is constant. The compactness of $X$ then implies that there exist $N$ and $R$ as above such that for all $x = (x_1, x_2, \ldots) \in X$, the value $f(x)$ is determined by the restrictions of $x_1$, $x_2$, \ldots, $x_N$ to $R$. In other words, there is a mapping $\rho \colon (6^R)^N \to k$ such that for all $x = (x_1, x_2, \ldots) \in X$,
 		\begin{equation}\label{eq:rho_box}
 			f(x) \,=\, \rho\left((x_i(\delta))_{1 \leq i \leq N, \delta \in R}\right).
 		\end{equation}
		We can then use \eqref{eq:rho_box} to define a continuous $\Pat$-avoiding $k$-coloring $f' \colon X_{\leq N} \to k$ of $X_{\leq N} \defeq \prod_{i = 1}^N X_i$.
		
		Now we claim that the conclusion of Lemma~\ref{lemma:F_loc_1} holds with \[F \,\defeq\, \left(\set{\mathbf{1}, \gamma_1, \ldots, \gamma_N} \cup R\right)^{4}.\] Take any $n \geq 2$ and a finite set $D \supseteq F$. For $1 \leq i \leq N$, let \[F_i^\ast \,\defeq\, \set{\mathbf{1}, \gamma_i}^{\log^\ast n + 2}.\] Then, by Lemma~\ref{lemma:F_loc_col_1}, the Schreier graph $G(X_{D,n}, \set{\gamma_i})$ admits an $F_i^\ast$-local proper $6$-coloring $f_i \colon X_{D, n} \to 6$. Define a $\G$-equivariant map $\pi_i \colon X_{D,n} \to X_i$ by $\pi_i(x)(\gamma) \defeq f_i(\gamma \cdot x)$ for all $x \in X_{D,n}$ and $\gamma \in \G$. Then
 		\[
 			\pi \colon X_{D,n} \to X_{\leq N} \colon x \mapsto (\pi_1(x), \ldots, \pi_N(x))
 		\]
 		is a $\G$-equivariant map from $X_{D,n}$ to $X_{\leq N}$. Thus, $f' \circ \pi \colon X_{D,n} \to k$ is a $\Pat$-avoiding $k$-coloring of $X_{D,n}$. Furthermore, to determine $(f' \circ \pi)(x)$, we only need to know $f_i(\delta \cdot x)$ for all $1 \leq i \leq N$ and $\delta \in R$, so this coloring is $(F_1^\ast \cup \ldots \cup F_N^\ast)R$-local. And now we are done since $F^\ast \defeq F^{\log^\ast n} \supseteq (F_1^\ast \cup \ldots \cup F_N^\ast)R$.
 	\end{scproof}
 
 	\subsection{Reduction to finite graphs}\label{subsec:tiles}
	
	For this subsection, we fix a finite subset $S \subset \G$. For each finite set $D \subset \G$ with $S \cup S^{-1} \cup \set{\mathbf{1}} \subseteq D$, we define a finite $S$-labeled graph $H_{D,n}$ as follows. For sets $A$ and $B$, let $\mathrm{Inj}(A,B)$ denote the set of all injective mappings from $A$ to $B$. The vertex set of $H_{D,n}$ is $V(H_{D,n}) \defeq \mathrm{Inj}(D, n)$. If $q$, $q' \in \mathrm{Inj}(D,n)$, we put an edge labeled $\sigma \in S \cup S^{-1}$ going from $q$ to $q'$ if and only if the following holds:
	\begin{equation}\label{eq:compat}
		\forall \delta,\, \delta' \in D, \qquad \left(\delta = \delta'\sigma \quad\Longrightarrow\quad q(\delta) = q'(\delta') \right).
	\end{equation}
	If \eqref{eq:compat} holds, we say that $q$ and $q'$ are \emphd{$\sigma$-compatible}. If $q$ and $q'$ are $\sigma$-compatible, then, in particular, $q'(\mathbf{1}) = q(\sigma)$. Since $q$ is injective, this implies that $q' \neq q$ and also that $q$ and $q'$ are not $\tau$-compatible for any $\tau \neq \sigma$, so the edge from $q$ to $q'$ in $H_{D,n}$ receives a unique label.
	
	\begin{lemma}\label{lemma:hom_exists}
		Let $D \subset \G$ be a finite set with $S \cup S^{-1} \cup \set{\mathbf{1}} \subseteq D$ and let $n \geq |D|^2$ be an integer. Then for every free zero-dimensional Polish $\G$-space $X$, there is a continuous homomorphism $G(X, S) \to H_{D,n}$.
	\end{lemma}
	\begin{scproof}
		The Schreier graph $G(X, D^{-1}D)$ has maximum degree at most $|D^{-1}D| - 1 \leq n - 1$ \ep{we are subtracting $1$ since $\mathbf{1} \in D^{-1}D$ does not count toward the degree}, so, by Lemma~\ref{lemma:coloring}, $G(X, D^{-1}D)$ has a continuous proper $n$-coloring $f \colon X \to n$. For each $x \in X$, let $q_x \colon D \to n$ be given by $q_x(\delta) \defeq f(\delta \cdot x)$ for all $\delta \in D$. By the choice of $f$, $q_x \in \mathrm{Inj}(D, n)$. Furthermore, it is clear that for any $\sigma \in S$, $q_x$ and $q_{\sigma \cdot x}$ are $\sigma$-compatible. Therefore, $x \mapsto q_x$ is a continuous homomorphism from $G(X,S)$ to $H_{D,n}$, as desired.
	\end{scproof}

	\begin{lemma}[\textls{Colorings of $H_{D,n}$}]\label{lemma:H_col}
		Let $\Pat$ be a finite set of $S$-connected $k$-patterns such that every free zero-dimensional Polish $\G$-space admits a continuous $\Pat$-avoiding $k$-coloring. Then there is a finite set $F \subset \G$ containing $S \cup S^{-1} \cup \set{\mathbf{1}}$ with the following property:
		
		\smallskip
		
		Let $n \geq 2$ and let $D \subset \G$ be a finite set. Set $F^\ast \defeq F^{\log^\ast n}$ and suppose that $D \supseteq F^\ast$. If $n \geq 2|D|$, then $H_{D,n}$ admits a $\Pat$-avoiding $k$-coloring.
	\end{lemma}
	\begin{scproof}
		Without loss of generality, we may assume that $\mathbf{1} \in \dom(p)$ for all $p \in \Pat$. Since each $p \in \Pat$ is $S$-connected, we can define $\Delta_p$ to be the diameter of $\dom(p)$ in $G(\G, S)$, i.e., the maximum length of a shortest path in $G(\dom(p), S)$ between two elements of $\dom(p)$. Set $\Delta \defeq \max_p \Delta_p$. Let $F \subset \G$ be given by Lemma~\ref{lemma:F_loc_1} applied to $\Pat$ and set
		\begin{equation}\label{eq:F}
			F \,\defeq\, (F_0 \cup \set{\mathbf{1}}) \left(S\cup S^{-1} \cup \set{\mathbf{1}}\right)^\Delta.
		\end{equation}
		Take any $n \geq 2$ and suppose that $D \supseteq F^\ast \defeq F^{\log^\ast n}$. Let  $F_0^\ast \defeq F_0^{\log^\ast n}$. Then, by Lemma~\ref{lemma:F_loc_1}, $X_{D,n}$ has an $F_0^\ast$-local $\Pat$-avoiding $k$-coloring $f \colon X_{D,n} \to k$, i.e., there is a map $\rho \colon n^{F_0^\ast} \to k$ such that for each $x \in X_{D,n}$,
		\begin{equation}\label{eq:rho_on_X}
			f(x) \,=\, \rho\left((x(\delta))_{\delta \in F_0^\ast}\right).
		\end{equation}
		We can simply use formula \eqref{eq:rho_on_X} to define a $k$-coloring $g$ of $H_{D,n}$; that is, for all $q \in \mathrm{Inj}(D,n)$, we let
		\begin{equation}\label{eq:rho_on_H}
			g(q) \,\defeq\, \rho\left((q(\delta))_{\delta \in F_0^\ast}\right).
		\end{equation}
		\ep{Here we are using that $F_0^\ast \subseteq F^\ast \subseteq D$.} We claim that $g$ is $\Pat$-avoiding, as desired.
		
		\begin{claim*}
			If $q \in \mathrm{Inj}(D,n)$, then there is a point $x \in X_{D,n}$ such that $x(\delta) = q(\delta)$ for all $\delta \in D$.
		\end{claim*}
		\begin{stepproof}
			The maximum degree of the Cayley graph $G(\G, D)$ is at most $|D \cup D^{-1}| - 1 \leq 2|D| - 1$ \ep{we are subtracting $1$ since $\mathbf{1} \in D$ does not count toward the degree}. Since $n \geq 2|D|$, we conclude that $q \colon D \to n$ can be extended to a proper $n$-coloring $x \colon \G \to n$ of $G(\G, D)$ greedily.
		\end{stepproof}
		
		Suppose that there is a pattern $p \in \Pat$ that occurs in $g$. This means that there is a homomorphism $\phi \colon \dom(p) \to \mathrm{Inj}(D,n)$ from $G(\dom(p), S)$ to $H_{D,n}$ such that $g(\phi(\gamma)) = p(\gamma)$ for all $\gamma \in \G$. By the above claim, there is a point $x \in X_{D,n}$ such that $x(\delta) = \phi(\mathbf{1})(\delta)$ for all $\delta \in D$. Since $\dom(p) \subseteq (S \cup S^{-1} \cup \set{\mathbf{1}})^\Delta$, equations \eqref{eq:F}, \eqref{eq:rho_on_X}, and \eqref{eq:rho_on_H} and the definition of $H_{D,n}$ yield $f(\gamma \cdot x) = g(\phi(\gamma))$ for all $\gamma \in \dom(p)$, so $p$ occurs in $f$, which is a contradiction.
	\end{scproof}

	Theorem~\ref{theo:compact} follows immediately from Lemmas~\ref{lemma:hom_exists} and \ref{lemma:H_col}. Fix an arbitrary increasing sequence $S \cup S^{-1} \cup \set{\mathbf{1}} \subseteq F_0 \subset F_1 \subset \ldots$ of finite subsets of $\G$ such that $\bigcup_{i = 0}^\infty F_i = \G$. Let $n_i \geq 2$ be any integer with
	\[
	n_i \,\geq\, |F_i|^{2\log^\ast n_i}.
	\]
	Set $D_i \defeq F_i^{\log^\ast n_i}$ and let $\mathscr{H} \defeq \set{H_{D_i, n_i}}_{i=0}^\infty$. Then Theorem~\ref{theo:compact} holds for this $\mathscr{H}$:
	
	\begin{theocopy}{theo:compact}
		Let $\Pat$ be a finite set of $S$-connected $k$-patterns. The following statements are equivalent:
		\begin{itemize}
			\item[\ref{item:all}] Every free zero-dimensional Polish $\G$-space admits a continuous $\Pat$-avoiding $k$-coloring.
			
			\item[\ref{item:compact1}] There is a graph in $\mathscr{H}$ that admits a $\Pat$-avoiding $k$-coloring.
			
			\item[\ref{item:compact_many}] All but finitely many graphs in $\mathscr{H}$ admit $\Pat$-avoiding $k$-colorings.
		\end{itemize}
	\end{theocopy}
	\begin{scproof}
		Implication \ref{item:compact_many} $\Longrightarrow$ \ref{item:compact1} is trivial, while \ref{item:compact1} $\Longrightarrow$ \ref{item:all} holds by Lemma~\ref{lemma:hom_exists} since $n_i \geq |D_i|^2$ for all $i$. Assuming \ref{item:all}, let $F \subset \G$ be given by Lemma~\ref{lemma:H_col} applied to $\Pat$. Then \ref{item:compact_many} holds since for all but finitely many $i$, we have $F_i \supseteq F$.
	\end{scproof}

	\subsection{\LOCAL algorithms}\label{subsec:LOCAL}
	
	In this subsection we prove Theorem~\ref{theo:LOCAL}:
	
	\begin{theocopy}{theo:LOCAL}
		Let $S \subset \G$ be a finite set and let $\Pat$ be a finite set of $S$-connected $k$-patterns. The following statements are equivalent: 
		\begin{itemize}
			\item[\ref{item:all}] Every free zero-dimensional Polish $\G$-space admits a continuous $\Pat$-avoiding $k$-coloring.
			\item[\ref{item:LOCAL}] There is a deterministic distributed algorithm in the \LOCAL model that, given an $n$-vertex $S$-labeled subgraph $G$ of $G(\G, S)$, in $O(\log^\ast n)$ rounds outputs a $\Pat$-avoiding $k$-coloring of $G$.
		\end{itemize}
	\end{theocopy}
		Implication \ref{item:LOCAL} $\Longrightarrow$ \ref{item:all} is a special case of \cite[Theorem~2.13]{BerDist}, so we only need to prove \ref{item:all} $\Longrightarrow$ \ref{item:LOCAL}. Before we proceed, let us record the following classical result, dating back to Goldberg, Plotkin, and Shannon \cite{GPS}, which can be seen as a distributed computing analog of Lemma~\ref{lemma:coloring}:
		
		\begin{theo}[{\cite[Corollary 3.15]{BE}}]\label{theo:d+1}
 		    There is a deterministic \LOCAL algorithm that computes a proper $(d + 1)$-coloring of an $n$-vertex graph $G$ of maximum degree $d$ in $\log^\ast n + O(d^2)$ rounds.
 	    \end{theo}
		
		Assume \ref{item:all} and let $F \subset \G$ be given by Lemma~\ref{lemma:H_col} applied to $\Pat$. Take $m$ so large that
		\begin{equation}\label{eq:m_lower}
			m \,>\, |F|^{3\log^\ast m}.
		\end{equation}
		Set $D \defeq F^{\log^\ast m}$. By Lemma~\ref{lemma:H_col}, the graph $H_{D,m}$ admits a $\Pat$-avoiding $k$-coloring $h \colon \mathrm{Inj}(D,m) \to k$. Thus, to prove \ref{item:LOCAL}, it suffices to show that there is a deterministic \LOCAL algorithm that, given an $n$-vertex $S$-labeled subgraph $G$ of $G(\G, S)$, in $O(\log^\ast n)$ rounds outputs a homomorphism $G \to H_{D,m}$ \ep{since composing such a homomorphism with $h$ requires no additional rounds of communication}.
		
		Our algorithm is supposed to output a homomorphism $G \to H_{D,m}$. In other words, each vertex $x \in V(G)$ has to compute an injective mapping $q_x \colon D \to m$ so that if $x$ is joined to $y$ by an edge with label $\sigma$, then $q_x$ and $q_y$ are $\sigma$-compatible. It is tempting to employ the same strategy as in the proof of Lemma~\ref{lemma:hom_exists}, i.e., to first compute, using Theorem~\ref{theo:d+1}, a locally injective $m$-coloring of $G$ and then make each vertex $x$ collect the colors within some finite radius around $x$ in $G$. Unfortunately, this approach does not quite work, because the set $D x$ may not be a subset of $V(G)$. Furthermore, there may be some $y \in D x$ whose distance to $x$ in $G$ is much larger than in $G(\G, S)$, so $x$ cannot find out the color of $y$ within a small number of rounds. We circumvent this difficulty by computing a homomorphism $G \to H_{D,m}$ directly.
		
		Let us start by introducing some useful notation. Let $\lambda$ be the edge labeling on $G$. We say that pairs $(x,\delta)$, $(y,\delta') \in V(G) \times D$ are \emphd{one-step equivalent}, in symbols $(x,\delta) \sim_1 (y, \delta')$, if $x$ and $y$ are adjacent and $\delta = \delta' \lambda(x,y)$. Note that the relation $\sim_1$ is symmetric. The equivalence relation on $V(G) \times D$ generated 
		by $\sim_1$ is denoted by $\sim$. Explicitly, we have $(x,\delta) \sim (y, \delta')$ if and only if either $(x, \delta) = (y, \delta')$ or there exists a finite sequence $(z_1, \delta_1)$, \ldots, $(z_t, \delta_t)$ such that \[(x, \delta) \sim_1 (z_1, \delta_1) \sim_1 \cdots \sim_1 (z_t, \delta_t) \sim_1 (y, \delta').\]
		When $(x,\delta) \sim (y, \delta')$, we say that $(x,\delta)$ and $(y, \delta')$ are \emphd{equivalent}. Observe that a mapping $x \mapsto q_x$ is a homomorphism from $G$ to $H_{D,m}$ if and only if $q_x(\delta) = q_y(\delta')$ whenever $(x,\delta) \sim (y, \delta')$.
		
		Let us establish a few simple facts about the relation $\sim$.
		
		\begin{big_claim}\label{claim:one}
			The following statements are valid:
			\begin{enumerate}[label={\ep{\itshape\alph*}}]
				\item\label{item:y} For every $x$ and $\delta$, $\delta' \in D$, there is at most one $y \in V(G)$ such that $(x,\delta) \sim (y, \delta')$.
				
				\item\label{item:delta} For every $x$, $y \in V(G)$ and $\delta \in D$, there is at most one $\delta' \in D$ such that $(x,\delta) \sim (y, \delta')$.
			\end{enumerate}
		\end{big_claim}
		\begin{scproof}
			Recall that $G$ is a subgraph of the Cayley graph $G(\G, S)$. Therefore, $x$ and $y$ are elements of the group $\G$, and if $(x,\delta) \sim (y, \delta')$, then we can write $\delta x = \delta' y$, so $y = (\delta')^{-1}\delta x$ and $\delta' = \delta x y^{-1}$.
		\end{scproof}
	
			
	
		For $x \in V(G)$, let $[x] \defeq \set{y \in V(G) \,:\, \text{$(x, \delta) \sim (y, \delta')$ for some $\delta$, $\delta ' \in D$}}$. Note that the relation ``$y \in [x]$'' is reflexive and symmetric, but not necessarily transitive.
	
		\begin{big_claim}\label{claim:close}
			For every $x \in V(G)$ and $y \in [x]$, the graph distance between $x$ and $y$ in $G$ is at most $|D|$.
		\end{big_claim}
		\begin{scproof}
			If $y = x$, then we are done. Otherwise, there is a sequence $(z_1, \delta_1)$, \ldots, $(z_t, \delta_t)$ such that
			\begin{equation}\label{eq:chain}
				(x, \delta) \sim_1 (z_1, \delta_1) \sim_1 \cdots \sim_1 (z_t, \delta_t) \sim_1 (y, \delta').
			\end{equation}
			By minimizing $t$, we may assume that the pairs $(x, \delta)$, $(z_1, \delta_1)$, \ldots, $(z_t, \delta_t)$, $(y, \delta')$ are pairwise distinct. By Claim~\ref{claim:one}\ref{item:y}, this implies that the elements $\delta$, $\delta_1$, \ldots, $\delta_t$, $\delta'$ are also pairwise distinct. Therefore, $t + 2 \leq |D|$. From \eqref{eq:chain}, we see that the distance between $x$ and $y$ is at most $t+1 < |D|$.
		\end{scproof}
	
		Let $G'$ denote the graph with the same vertex set as $G$ in which two distinct vertices $x$, $y$ are adjacent if and only if there is $z \in V(G)$ such that $z \in [x]$ and $y \in [z]$ \ep{this includes the case when $z = x$ and $y \in [x]$}. By Claim~\ref{claim:one}\ref{item:y}, $|[x]| \leq |D|^2$ for all $x \in V(G)$, so the maximum degree of $G'$ is at most
		\[
			N \,\defeq\, |D|^4 \,=\, O(1).
		\]
		\ep{Here and in what follows, the asymptotic notation is with respect to $n \to \infty$.} By Claim~\ref{claim:close}, a single communication round in the \LOCAL model on $G'$ can be simulated by $2|D| = O(1)$ rounds on $G$. Hence, by Theorem~\ref{theo:d+1}, we can compute a proper $(N + 1)$-coloring $\phi \colon V(G) \to (N + 1)$ of $G'$ in $O(\log^\ast n)$ rounds. For $0 \leq i \leq N$, let $X_i \defeq \phi^{-1}(i)$. 
		
		We shall compute the desired homomorphism from $G$ to $H_{D,m}$ in $N + 1$ stages indexed by $0$, $1$, \ldots, $N$. At the start of stage $i$, each vertex $x$ will have already computed the values $q_x(\delta)$ for some subset of $\delta \in D$, subject to the following requirement:
		\begin{quote}
			\emph{If $(x, \delta) \sim (y, \delta')$, then $q_x(\delta) = q_y(\delta')$ whenever at least one of $q_x(\delta)$ and $q_y(\delta')$ is defined.}
		\end{quote}
		During stage $i$, we have to compute $q_x(\delta)$ for all $x \in X_i$ and $\delta \in D$. To this end, each vertex $x \in X_i$ considers the elements $\delta \in D$ one by one and performs the following procedure for each of them. If $q_x(\delta)$ is already defined, then there is nothing to do. Otherwise, by Claim~\ref{claim:close}, in $|D|$ rounds $x$ can determine the following set:
		\[
			B \, \defeq\, \set{q_y(\epsilon) \,:\, \text{$y \in [x]$, $\epsilon \in D$, and $q_y(\epsilon)$ is defined}}.
		\]
		Since $|[x]| \leq |D|^2$, we have $|B| \leq |D|^3 < m$ by \eqref{eq:m_lower}, so $x$ can pick a color $\alpha < m$ that is not in $B$ and set $q_x(\delta) \defeq \alpha$. Then in $|D|$ rounds $x$ can notify each $y \in [x]$ so that if $(x,\delta) \sim (y, \delta')$, then $y$ sets $q_y(\delta') \defeq \alpha$ \ep{such $\delta'$ is unique by Claim~\ref{claim:one}\ref{item:delta}}. By the choice of $\alpha$, the mappings $q_y \colon D \rightharpoonup m$ for all $y \in [x]$ remain injective after this procedure. Notice also that since the set $X_i$ is $G'$-independent, all the elements of $X_i$ can run this procedure in parallel without creating any conflicts.
		
		After $(N+1)$ stages, we will have computed a homomorphism $G \to H_{D, m}$. Note that each stage takes $O(1)$ rounds, and there are $O(1)$ stages, so the total required number of communication rounds is
		\[
			\underbrace{O(\log^\ast n)}_{\text{computing $\phi$}} \,+\, O(1) \,=\, O(\log^\ast n),
		\]
		and the proof is complete.

    \appendix
    
    \section*{Appendix: Proof of Lemma~\ref{lemma:F_loc_col_1}}\label{app}
    
    Here we give a proof of Lemma~\ref{lemma:F_loc_col_1}:
    
    \begin{lemmacopy}{lemma:F_loc_col_1}
 		Let $\gamma \in \G \setminus \set{\mathbf{1}}$ and let $D \subset \G$ be a finite set with $\gamma \in D$. Take $n \geq 2$ and let $F^\ast \defeq \set{\mathbf{1}, \gamma}^{\log^\ast n + 2}$. Then the Schreier graph $G(X_{D,n}, \set{\gamma})$ admits an $F^\ast$-local proper $6$-coloring.
 	\end{lemmacopy}
    
 	    As $\gamma \in D$, we have $X_{D, n} \subseteq X_{\set{\gamma}, n}$, so it is enough to consider the case when $D = \set{\gamma}$. For brevity, let
 	    \[
 	        X_n \,\defeq\, X_{\set{\gamma}, n} \quad \text{and} \quad G_n \,\defeq\, G(X_{\set{\gamma}, n}, \set{\gamma}).
 	  \]
 	    If $n \leq 6$, then the mapping $X_n \to 6 \colon x \mapsto x(\mathbf{1})$ is as desired \ep{it is a $\set{\mathbf{1}}$-local proper $6$-coloring of $G_n$}. Thus, we may assume that $n > 6$. For every $n \geq 6$, we define \[f(n) \,\defeq\, 2 \lceil \log_2 n\rceil\] and observe that $f(6) = 6$ and $n > f(n) \geq 6$ for $n > 6$.
 	    
 	    The heart of the construction is in the following claim:
 	      
 	    \begin{claim*}
			For every $n > 6$, the graph $G_n$ admits a $\set{\mathbf{1}, \gamma}$-local proper $f(n)$-coloring.
		\end{claim*}
		\begin{scproof}
			For an integer $a$ with $0 \leq a \leq n-1$, let $a_0a_1\ldots a_{\lceil \log_2 n \rceil-1}$ denote the binary expansion of $a$. For each $x \in X_n$, let $i(x)$ be the smallest index $i$ such that $x(\mathbf{1})_i \neq x(\gamma)_i$ \ep{such an index exists since $x$ is a proper coloring of $G_n$} and set $d(x) \defeq x(\mathbf{1})_{i(x)}$ and
			$c(x) \defeq 2i(x) + d(x)$. By construction, the function $c$ is $\set{\mathbf{1}, \gamma}$-local and $0 \leq c(x) < 2\lceil \log_2 n \rceil =  f(n)$. It remains to verify that $c$ is a proper coloring of $G_n$. Suppose, toward a contradiction, that $c(x) = c(\gamma \cdot x)$ for some $x \in X_n$. Then $i(x) = i(\gamma \cdot x)$, because both these quantities are equal to $\lfloor c(x)/2\rfloor$. Letting $i \defeq i(x)$, we see that $d(x) = x(\mathbf{1})_i \neq x(\gamma)_i = d(\gamma \cdot x)$ by the definition of $i(x)$. But then $c(x) \neq c(\gamma \cdot x) \pmod 2$; a contradiction.
		\end{scproof}
		
		For $n > 6$, let $f^\ast(n)$ be the minimum $k$ such that $f^{(k)}(n) = 6$, where
 	    $
 	        f^{(k)} \defeq f \circ \cdots \circ f
 	      $ \ep{$k$ times}.
		It is routine to check that $f^\ast(n) \leq \log^\ast n + 2$ for all $n > 6$. By iterating the above claim $f^\ast(n)$ times, we obtain a sequence of $\G$-equivariant maps
		\[
		    X_n \,\xrightarrow{\quad \pi_1 \quad}\, X_{f(n)} \,\xrightarrow{\quad \pi_2 \quad }\, X_{f^{(2)}(n)} \,\xrightarrow{\quad \pi_3\quad }\, \cdots \,\xrightarrow{\ \ \pi_{f^\ast(n)}\ \ }\, X_{6},
		\]
		where for each $i$, the mapping $x \mapsto \pi_i(x)(\mathbf{1})$ is $\set{\mathbf{1}, \gamma}$-local. It remains to set $\pi \defeq \pi_{f^\ast(n)} \circ \cdots \circ \pi_1$ and observe that the mapping $x \mapsto \pi(x) (\mathbf{1})$ is an $F^\ast$-local proper $6$-coloring of $G_n$, as desired.

	\printbibliography
	
\end{document}